%% file: ArXiv.tex
\documentclass[final,leqno]{siamltex}
\usepackage[latin1]{inputenc}
\usepackage[english]{babel}
\usepackage{graphicx}
\usepackage{amsmath,amsfonts,amssymb,amsopn,amscd}
\usepackage{bbm,wasysym,stmaryrd}
\usepackage{subfigure}
\usepackage{color}
\usepackage{hyperref}
\usepackage[normalem]{ulem}

\newcommand{\R}{\mathbb{R}}
\newcommand{\C}{\mathbb{C}}
\newcommand{\J}{\mathcal{J}}
\newcommand{\I}{\mathcal{I}}
\newcommand{\M}{\mathcal{M}}
\newcommand{\U}{\mathcal{U}}
\newcommand{\erf}{\textrm{erf}}
\renewcommand{\S}{\mathbb{S}}

\newcommand{\bqq}{\begin{equation}}
\newcommand{\eqq}{\end{equation}}
\newcommand{\bqs}{\begin{equation*}}
\newcommand{\eqs}{\end{equation*}}

\newcommand{\N}{\mathbb{N}}

\newcommand{\E}{\mathcal{E}}
\newcommand{\G}{\mathcal{G}}

\newcommand{\cl}{\mathcal{L}}
\newcommand{\V}{\mathcal{V}}

\newcommand{\mx}{\mathbf{x}}
\newcommand{\my}{\mathbf{y}}
\newcommand{\mz}{\mathbf{z}}
\newcommand{\me}{\mathbf{e}}

\newcommand{\mG}{\mathbf{G}}

\newcommand{\vp}{\varphi}

\usepackage{graphicx}
\usepackage{subfigure}
\graphicspath{{./Figures/}}

\title{Pulsatile localized dynamics in delayed neural-field equations in arbitrary dimension}

\author{Gr\'egory Faye\thanks{University of Minnesota,
School of Mathematics,
206 Church Street S.E.,
Minneapolis, MN 55455 US  (\href{mailto:gfaye@umn.edu}{gfaye@umn.edu}). This author was supported by an NSF grant (DMS-1311414).} \and  Jonathan Touboul\thanks{Mathematical Neuroscience Lab, CIRB-Coll\`ege de France, 11 place Marcelin Berthelot, 75005 Paris (\href{mailto:jonathan.touboul@college-de-france.fr}{jonathan.touboul@college-de-france.fr}) \& INRIA Paris-Rocquencourt, Mycenae Team. } }

\begin{document}

\maketitle

\begin{abstract}
Neural field equations are integro-differential systems describing the macroscopic activity of spatially extended pieces of cortex. In such cortical assemblies, the propagation of information and the transmission machinery induce communication delays, due to the transport of information (propagation delays) and to the synaptic machinery (constant delays). We investigate the role of these delays on the formation of structured spatiotemporal patterns for these systems in arbitrary dimensions. We focus on localized activity, either induced by the presence of a localized stimulus (pulses) or by transitions between two levels of activity (fronts). Linear stability analysis allows to reveal the existence of Hopf bifurcation curves induced by the delays, along different modes that may be symmetric or asymmetric. We show that instabilities strongly depend on the dimension, and in particular may exhibit transversal instabilities along invariant directions. These instabilities yield pulsatile localized activity, and depending on the symmetry of the destabilized modes, either produce spatiotemporal \emph{breathing} or \emph{sloshing} patterns.
\end{abstract}

\begin{keywords} 
Neural fields equations, integro-differential equations, delays, Hopf bifurcation, breathers.
\end{keywords}

\begin{AMS} 92C20; 35R09; 34K13.   \end{AMS}

\pagestyle{myheadings}
\thispagestyle{plain}

\section{Introduction}
\input Introduction.tex
\section{Neural Fields Formalism with Delays}\label{sec:NFE}
\input Amari

\section{Pulses: Existence, Stability and Delay-Induced Breathers and Sloshers}\label{sec:Pulses}
\input Pulse

\section{Breathing fronts in general neural fields}\label{sec:Fronts}
\input Front

\section{Discussion}
\input Conclusion

\section*{Aknowledgements} The authors would like to warmly thank Richard Bailleul for several insightful discussions and help at the early stage of this work, in the frame of his undergraduate internship at the Mathematical Neuroscience Lab.

\bibliographystyle{plain}


\end{document}

%% file: Introduction.tex
Electrode recordings and imaging studies have revealed that the cortex can support a variety of spatiotemporal dynamics including standing waves~\cite{BFC,SNC}, traveling pulses~\cite{BFC,HTYMLSW,WHZ} and spiral waves~\cite{HTYMLSW,HXLTGW,WHZ}. These patterns of activity are visible from macroscopic measurements related to the firing rate of neurons at a specific place, and is therefore the macroscopic result of the activity of individual cells that compose the cortical tissue. The interaction between nerve cells takes a specific delay: action potentials are delivered to the interconnected neurons after a time depending on the distance between cells and on the time taken by the synapse to transmit the information. These delays,  on the order of milliseconds, are due to several phenomena including the finite-velocity propagation of action potentials, dendritic and synaptic processing as well as spike-generation dynamics. The question we may address in the present manuscript is the role of these delays in shaping the macroscopic dynamics of large neuronal assemblies.

The first step in this journey is to describe the effective macroscopic of large neuronal networks. In that domain, a successful approach consists in considering that the activity of large interconnected networks satisfies a \emph{Neural Field Equation} (NFE). These models introduced in the end of last century~\cite{amari:77,wilson-cowan:73} describe the macroscopic activity of the network through an averaged variable, the mean firing rate, of a local neuronal population. This firing rate evolves depending on the input received by the population, which is the sum of an external current term (representing the external stimulations for instance) and the activity of cells at other locations on the cortex. Thus descriptions of spatially extended populations of neurons give rise to integro-differential evolution equations. These heuristic descriptions provide a foundation for studying large populations of cells with complex dynamics, and can be rigorously seen as accurate descriptions of the underlying microscopic dynamics of noisy spatially extended networks as shown in~\cite{T,T1,T2}. These neural fields equations were indeed very successful in reproducing a number of phenomena as reviewed in the comprehensive article~\cite{bressloff:12}, including in particular visual hallucination patterns~\cite{ermentrout-cowan:79,bressloff-cowan-etal:01}, binocular rivalry~\cite{bressloff2012neural}, working memory~\cite{laing-troy-etal:02}. From this wide literature, it became clear that, as is the case of the widely studied nonlinear reaction diffusion systems of PDEs (with which the neural field equation has several common points~\cite{bressloff2012invasion}), the class of neural field equation supports spatio-temporal pattern formation and may exhibit a wide repertoire of phenomena, such as waves, fronts, pulses and localized activity~\cite{bressloff:12}. 

The role of delays has been shown central in the activity of neural fields, in relationship with the emergence of traveling waves or standing waves, and how these delays may affect the stability of homogeneous solutions (see~\cite{veltz2013interplay} and references therein). All these studies point to a specific role of delays in shaping the macroscopic activity of neural fields, and tend to induce oscillations, as characterized by the exhibition of delay-induced instabilities. However, the role of delays in the stability of stimulus-locked pulses or localized fronts has not been addressed, and is of central importance in order to understand the role of delays in the activity of neural fields in the presence of localized input, or in the stability of fronts. 

In the absence of delays, localized stimulus generates a super-threshold response whose stability was characterized in the seminal works of Amari~\cite{amari:77}. The dynamics of cells is actually more complex, and particularly neuronal adaptation can yield destabilization of the pulse. This was demonstrated recently investigated~\cite{F,folias-bressloff:04} in the Pinto and Ermentrout model of neural field with adaptation~\cite{pinto-ermentrout:01,pinto-ermentrout:01b}. In detail, stationary pulses may undergo a Hopf bifurcation as the input amplitude is decreased, leading to stable periodic oscillations around the localized pulse. These oscillations take the form of \emph{breathers} (periodically expanding and contracting solutions around the stationary pulse) or \emph{sloshers} (periodic solutions asymmetrically varying around the stationary pulse). 

Here, we will be interested in the influence of delays on the emergence of such solutions, and we will extend these works to the case of fronts. In order to separate the effects due to linear adaptation and delays, we will work here in the classical Amari framework~\cite{amari:77}, in the presence of constant and space-dependent propagation delays. Phenomena arising in these equations are highly dependent of the dimension of the neural field. In order to take this aspect into account, theoretical developments will be performed in arbitrary dimensions, and examples in dimension one (simplest example) and two (more relevant to the underlying biological problem). In all cases, we will show that delays can destabilize the stationary pulses through Hopf bifurcations, yielding, similarly to the case of Pinto-Ermentrout model with adaptation, to breathing or sloshing pulses. We extend these analyses in order to include the case of stationary fronts in general dimensions. In that case again, we will show that both constant and propagation delays shape the solution of the neural field, yielding localized periodic oscillations of the front location. Fronts are essentially one-dimensional patterns, yet instabilities may arise in transverse dimensions, phenomenon only arising in dimension larger than two. All these analytical developments are confirmed by extensive numerical analyses that show a very good agreement with the predicted transition points. 

The paper is organized as follows. In section~\ref{sec:NFE}, we introduce the framework of the study. Section~\ref{sec:Pulses} will deal in the existence and stability of localized activity (pulses), and section~\ref{sec:Fronts} will consider the case of stationary fronts and their stability. 

%% file: Amari.tex
The classical Neural Field Equation describes a continuum limit of the activity of neurons belonging to a spatially extended piece of cortex $\Gamma \subset \R^d$ through a macroscopic variable, the population-averaged firing rate $u(\mx,t)$ at location $\mx\in\Gamma$ and time $t\geq 0$. This average firing rate evolves depending of the input it receives. In the absence of input, it is postulated that the activity exponentially decays to zeros. Input received by the neural field are of two types:
(i) the external input at time $t$ received at location $\mx\in\Gamma$, denoted $I(\mx,t)$, accounting for the external stimuli received by the cortex, and (ii) the internal input from neurons within the network. This second type of input necessitates to describe how the average firing-rate at location $\mx\in \Gamma$ is affected by the presence of an averaged firing rate at location $\mx'$. This depends on the structure of the network, in particular on the distance between $\mx$ and $\mx'$ (inducing propagation delays), and on the type of connection (level of excitation or inhibition) between the two areas. Here, we will consider the classical synaptically coupled network. In this framework, probably the most widely used neural field model since the seminal works of Amari and Wilson and Cowan~\cite{amari:77,wilson-cowan:72}, we will consider that the type of interaction between these populations is described by a kernel $\J_d(\|\mx-\mx'\|)$ which only depends on the distance between the cells, and that the effective input received from neurons at location $\mx'$ is a non-decreasing function of the mean firing rate $f(u(\mx',t-\tau(\mx-\mx')))$. This framework therefore considers a single layer composed of both excitatory and inhibitory neurons, and gives rise to the following integro-differential equation: 
\begin{equation}\label{eq:NFE}
\frac{\partial u}{\partial t}(\mx,t) = -u(\mx,t) + \int_{\Gamma} \J_d(\|\mx-\my\|)f(u(\my,t-\tau( \mx-\my)) d\my + I(\mx,t)
\end{equation}
The quantity $f(u(\mx,t))$ represents the average activity, namely the pulse emission rate, of the neurons at location $\mx$ and time $t$, it is therefore a non-decreasing positive function. The simplest model one could think of, and probably the only model allowing closed-form derivation of the solutions to the NFE (see e.g.~\cite{bressloff:12,amari:77,F}) is the Heaviside function with threshold $\theta$ : 
\begin{equation}\label{eq:Heaviside}
f(u)=
\left\{ 
\begin{array}{l}
0~~,~~u \leq \theta \\ 
1~~,~~u > \theta
\end{array} 
\right.
\end{equation}

\begin{figure}[htbp]
		\centering
			\subfigure[Mexican Hat]{\includegraphics[width=.4\textwidth]{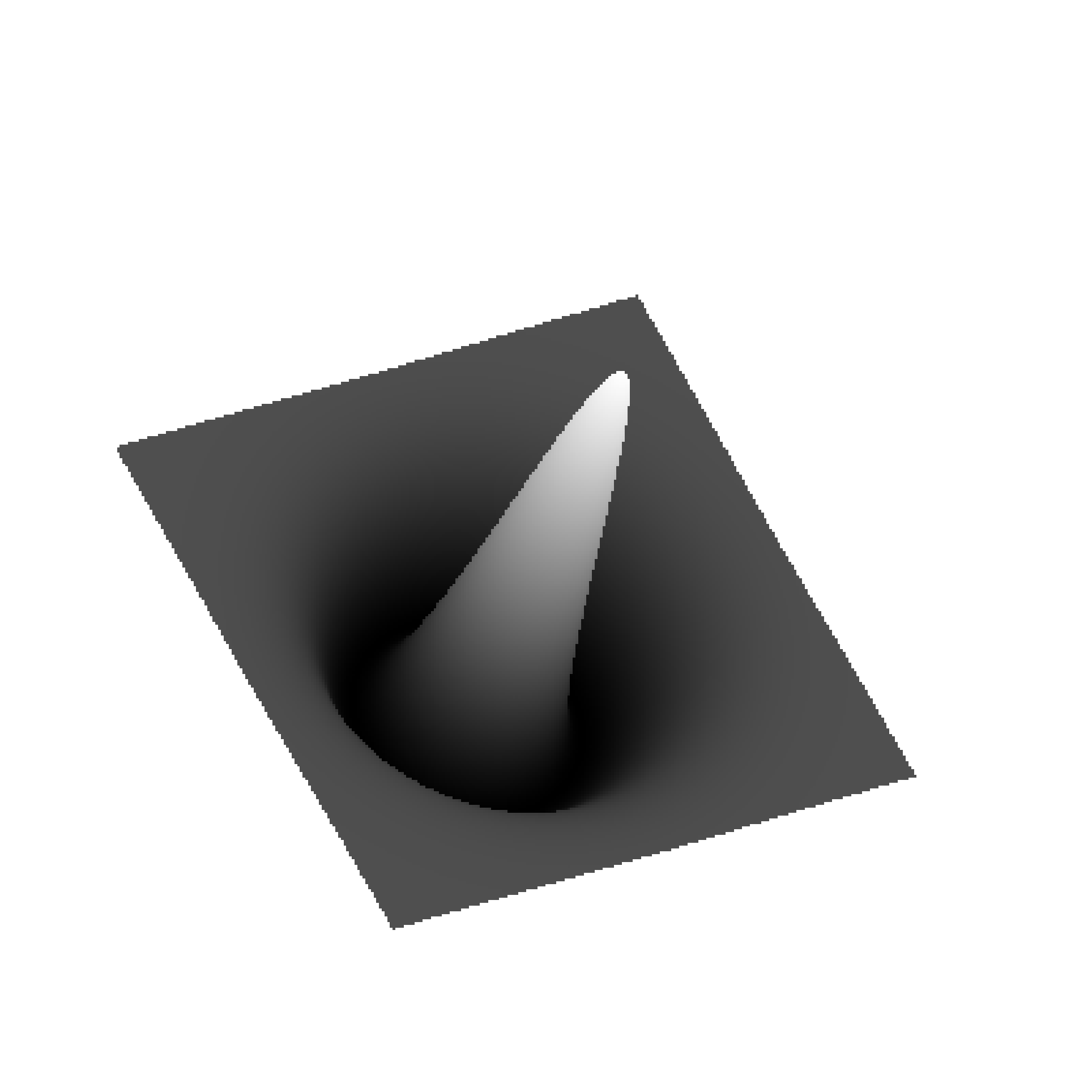}}
			\subfigure[Inverse Mexican Hat]{\includegraphics[width=.4\textwidth]{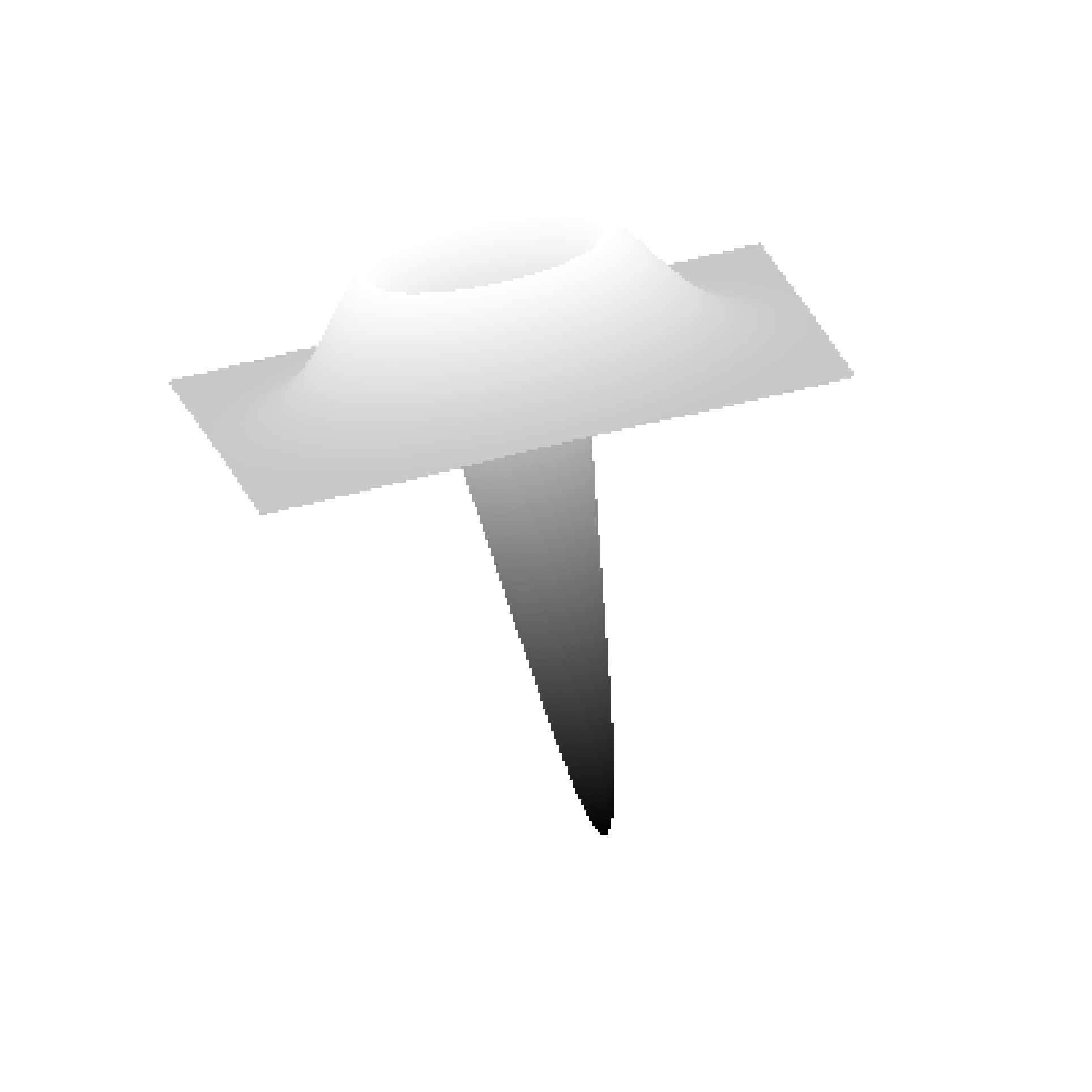}}
		\caption{Typical shape of the connectivity functions $\J_d$ in two dimensions ($\Gamma \subset \R^2$). Here $\J_d$ is a difference of Gaussians. (a) local excitation and lateral inhibition case (functional connectivity profile), (b) local inhibition and lateral excitation (anatomical profile).}
		\label{fig:MexHat}
\end{figure}

The connectivity function $\J_d$ represents the average intensity of connection between neurons at positions $\mx$ and $\my$ on the neural field $\Gamma$. It characterizes the averaged effect of spikes received by neurons at location $\mx$ incoming from neurons at location $\my$. The choice of this function therefore depends on the type of anatomical and functional connections between different neuronal areas. Typically, it has been shown that functional connectivity is globally excitatory at short spatial scales and inhibitory at larger distances. This macroscopic description of the connection motivated to consider functions that only depend on the distance on the cortex. Anatomical studies showed that excitatory axons are in average, longer than inhibitory cables, yielding to a typical \emph{inverse Mexican hat} connectivity. When $\Gamma$ is a feature space representing the coding attribute of neurons (e.g., orientation of a visual stimulus for the modeling of primary visual areas~\cite{bressloff-cowan-etal:02}), connections are locally excitatory and laterally inhibitory (Mexican hat connectivity, see Fig.~\ref{fig:MexHat}).
In the present study, we will see that the specific shape of the map $\J_d$ is not essential to obtain rigorous closed-form results on the existence and stability of particular types of solutions. For concreteness, we will consider that these maps are differences of Gaussians:
\[\J_d(\|\mx\|) = w_e \G_d(\|\mx\|,\sigma_e) - w_i \G_d(\|\mx\|,\sigma_i),\quad \forall \mx \in \R^d,\]
where $\G_d(\Vert \mx\Vert,\sigma)$ is the isotropic Gaussian function in dimension $d$ with variance $\sigma$:
 \[ \G_{d}(x,\sigma)= \frac{1}{(\sigma\sqrt{\pi})^d}\exp\left(-\frac{x^2}{\sigma^2}\right), \quad \forall x\in\R.\]

In equation~\eqref{eq:NFE}, the term $I(\mx,t)$ corresponds to the intensity of a stimulus reaching the neural field at location $\mx$ and time $t$. In section~\ref{sec:Pulses}, we will choose time-independent localized stimuli around the origin, and will opt for a Gaussian form: $I(\mx)=\I(\| \mx\|)=I_0e^{-\| \mx\|^2/\sigma^2}$. This will induce localized pulse activity. In section~\ref{sec:Fronts}, we will work with an independent stimulus of the form: 
\bqs
I(\mx)=\I(\me\cdot \mx)=I_0\left(1-\frac{1}{1+e^{-s \me\cdot\mx}} \right),
\eqs 
where $\me\in \mathbb{S}^{d-1}$ will be the direction of propagation for the fronts. Again, the particular shape of these two stimuli will not crucially impact the results. 

Eventually, the delay function is considered of the general form 
\[\tau (\mx) = \tau_D + \frac{\|\mx\|}{c}\]
such that it takes into account both constant transmission delays ($\tau_D$) and space-dependent delays ($c$) related to the transmission of spikes at finite speed along the axons.

The model is therefore complete as soon as the space $\Gamma$ is defined. Since we will investigate in the present manuscript localized states, we will consider that the piece of tissue considered, $\Gamma$, covers $\R^d$. Beyond the fact that under this assumption the problem is amenable to analysis, this appears as a reasonable assumption in the present context. Indeed, boundary effects may be negligible as long as the spatial extent of the phenomena investigated is relatively small compared to the size of the entire network, which is sensible when considering the presence of localized activity patterns or the local evolution of stationary fronts.

Now that these equations and assumptions have been done, we are in a position to analyze their solutions. We start by analyzing the presence of localized patterns of activity (pulses) driven by the presence of a local stimulus.

%% file: Pulse.tex
The basic mechanism described by Bressloff and Folias~\cite{folias-bressloff:04,F} for the generation of a stimulus induced breathing pulse in the Pinto-Ermentrout model consists in a Hopf instability of the stimulus-induced stationary pulse. We will show here that the same phenomenon arises, in arbitrary dimensions, related to the presence of delays. Beforehand, we will briefly discuss (in an overarching formalism) the existence of pulses in general dimensions.

\subsection{Existence of pulses in general dimensions}\label{sec:ExistPulse}
For one-dimensional equations, Amari provided the conditions for the existence of stationary pulses~\cite{amari:77}, and similar arguments were provided in two dimensions~\cite{coombes-owen:04}. These arguments readily extend to $d$-dimensional neural fields.

Let us consider equation~\eqref{eq:NFE} in dimension $d\geq 1$ and consider a rotationally invariant time-independent input $\I(\| \mx \|)$ centered at the origin. A stationary rotationally invariant pulse (simply called pulse in the sequel) is a time-independent profile $\mathcal{U}(\mx)$ solution of equation~\eqref{eq:NFE} satisfying, for some $a>0$ (called the pulse half-width):
\begin{equation}\label{eq:BumpDef}
\begin{cases}
	\mathcal{U}(\mx) > \theta, & \| \mx\| < a \\ 
	\mathcal{U}(\mx) < \theta, & \text{otherwise}\\
\end{cases}
\end{equation}
and decaying towards zero at infinity. 

Plugging these conditions in equation~\eqref{eq:NFE} yields the implicit relation: 
\begin{equation}\label{eq:bump}
\mathcal{U}(\mx)= \int_{\R^d} \J_d(\|\mx-\my\|) f\left (\mathcal{U}(\my) \right ) d\my+\I(\| x \|),
\end{equation}
which provides the profile of possible pulses parametrized by their half-width $a$:
\begin{equation}\label{eq:bumpHalfWidth}
\mathcal{U}(\mx)= \int_{B_d(0,a)} \J_d(\|\mx-\my\|) d\my+\I(\| \mx \|). 
\end{equation}

The question that remains to answer is whether there exists such solutions, or in other words if there exists such quantity $a$ so that equation~\eqref{eq:BumpDef} is satisfied. Since the functions are rotationally invariant, we have that for any $\mx$ such that $\| \mx \| = r$, the integral term in the righthand side is equal to:
\[\M(r,a) = \int_{B_d(0,a)} \J_d(\|\mx-\my\|)\,d\my = \int_0^a \int_{\mathbb{S}_{d-1}} \J_d( \|r- u \varphi\|) u^{d-1} du d_{d-1}\varphi, \]
where $B_d(0,a)$ is the ball of radius $a$ in $\R^d$ and $\mathbb{S}_{d-1}$ is the unit sphere. Therefore, existence of a pulse is ensured as soon as there exists $a>0$ satisfying the implicit relationship:
\begin{equation}\label{eq:BumpExistence}
	\theta = \M(a,a) + \I(a).
\end{equation}
This provides a condition for a bump to exist provided that it crosses the threshold $\theta$ exactly once. Moreover, finding a solution to this equation provides the exact shape of the solution of the pulse through equation~\eqref{eq:bump}. Eventually, this relationship provides the non-linear dependence of the half-width $a$ with respect to the other parameters. 

The existence of stimulus-dependent pulse solutions is of course dependent on the maximal amplitude of the input $I_0$. The pulse relationship can be easily written as a condition on this amplitude, and one obtains:
\[I_0 = e^{a^2/\sigma^2}(\theta-\M(a,a)).\]
Interestingly, we observe that the level of input necessary to induce the presence of a stationary bump, as well as its half-width, depends on the dimension of the neural field, and on the type of connectivity function chosen.

In dimension $1$, we obtain therefore the simple conditions for bumps to exist as:
\[\theta = \int_0^{2a} \J_1(x) dx + \I(a)\]
which in the case of our choice of functions further simplifies as:
\[\theta=w_e \erf\left(\frac{2a}{\sigma_e}\right)-w_i \erf\left(\frac{2a}{\sigma_i}\right) + I_0 e^{-a^2/\sigma^2}.\]

It is now convenient, in order to characterize more precisely the shape of the bump (and in order to expedite the analysis of the stability of the bump) to compute the gradient of this solution. By standard differential calculus methods, one obtains:
\begin{align}
	\label{eq:GradientBump}\nabla \U (\mx) &= \int_{B_d(0,a)} \J_d' (\| \mx-\my\|)\frac{\mx-\my}{\| \mx - \my\|}\,d\my + \nabla \I (\|\mx\|)\\
	\nonumber & = a^d \int_0^1 \int_{\S_{d-1}} \J_d' (\| \mx-a u \vp\vert)\frac{\mx-a u \vp}{\| x - a u \vp\|}\,u^{d-1}du d_{d-1}\vp + \I' (\| \mx\|)\frac{\mx}{\| \mx \|}.
\end{align}
Of course, as the gradient of a spherically symmetric function, this function is odd (i.e. $\nabla \U(-\mx)=-\nabla \U(\mx)$), and more generally, is equivariant under the action of the orthogonal group: for any orthogonal transformation (e.g., rotation) $A\in \mathbf{O}(d)$, we have:
\begin{align*}
	\nabla \U(A\mx) &= \int_{B_d(0,a)} \J_d' (\| A\mx-\my\|)\frac{A\mx-\my}{\| A\mx - \my\|}\,d\my + \nabla \I (\|A\mx\|)\\
	&=\int_{B_d(0,a)} \J_d' (\| \mx-\mz\|)\frac{A\mx-A\mz}{\| \mx - \mz\|}\,d\mz + \I'(\|\mx\|) \frac{A\mx}{\| \mx\|} \\
	&= A\nabla \U(\mx).
\end{align*}
This simple property will be useful in order to characterize the deformation modes of the bumps. In particular, it is important to note that the norm of this quantity is constant along spheres in $\R^d$, which will be useful for the analysis of the stability.

In dimension $d=1$, we obtain the rather classical and simple expression at $x=\pm a$
\[\U'(-a)= -\U'(a)=\J_1(0)-\J_1(2a) + \I'(a).\]

\begin{figure}[htbp]
		\centering
		\includegraphics[width=.9\textwidth]{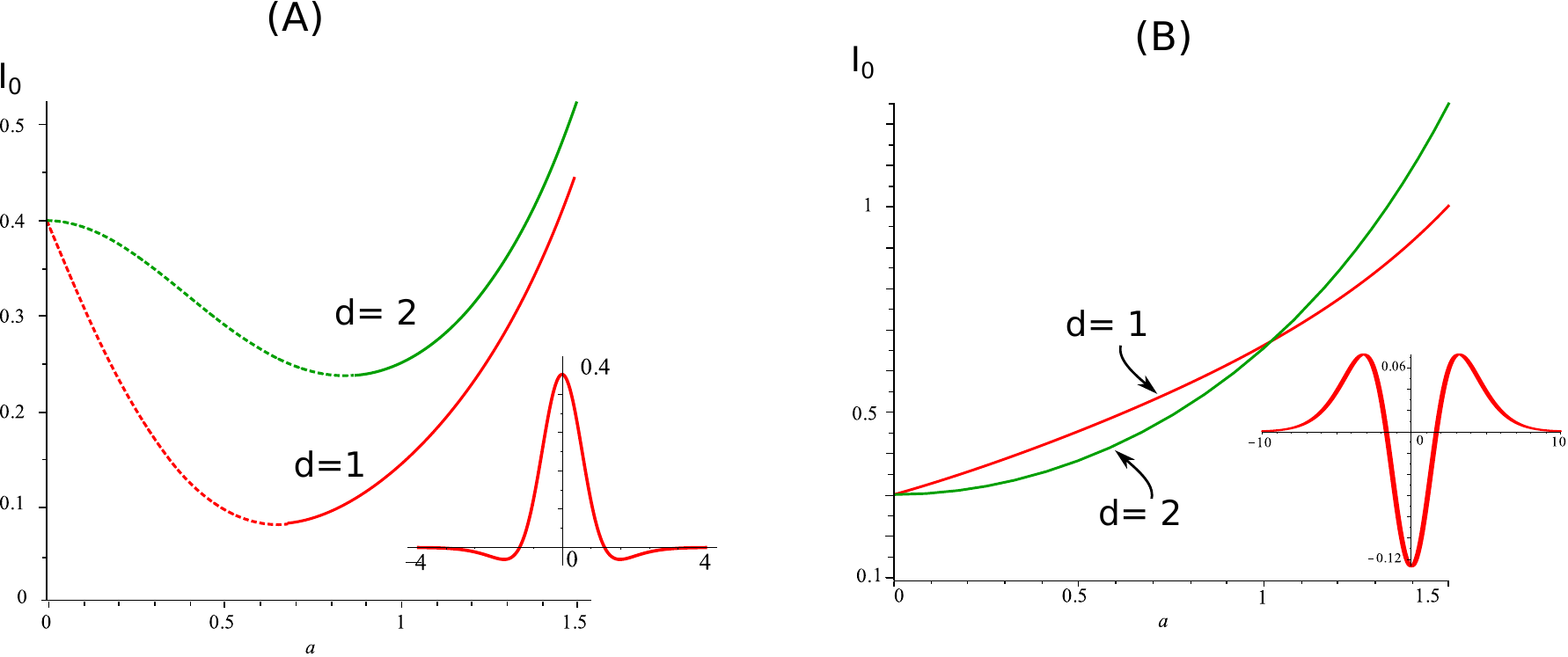}
		\caption{Dependence of the stationary pulse half-width a on the intensity of the input $I_0$ in the 1d (red line) and 2d (green line) neural fields. Dashed (solid) lines indicate the stability of the bump in the non-delayed case. (A) Mexican hat connectivity (see inset for the 1-d connectivity) $\theta=0.4$, $\sigma_e=1$, $\sigma_i=2$, $\sigma=1.2$, $w_e=0.4$ and $w_i=1$. (B) Inverse Mexican hat connectivity (see parameters in Fig.~\ref{fig:Bump1DSynapticDelayModelI}).}
		\label{fig:I0}
\end{figure}

Figure~\ref{fig:I0} represent the solutions to the bump equation for $d=1$ or $2$. We observe that the existence and size of the bump depends on the dimension and on the type of connectivity chosen.
In the case of the Mexican hat connectivity (Fig.~\ref{fig:I0}(A)), the value of the input amplitude $I_0$ forms a U-shaped curve as a function of the pulse size $a$. For $I_0$ smaller than a critical value (which is strictly smaller than the threshold $\theta$), no bump exist, and when $I_0$ exceeds this value, two bumps exist for $I_0<\theta$, and a unique bump exists for $I_0>\theta$. At $I_0=\theta$, the size of the smallest bump goes to zero. The existence of pulses for inputs strictly smaller than the threshold is precisely due to the local excitation of Mexican hat connectivities. In that case, some values of the input sufficient to trigger a pulse in dimension $1$ are too small to trigger a pulse in dimension $2$. In the case of the inverse Mexican hat connectivity (Fig.~\ref{fig:I0}(a)), this is no more the case: the value of the half-width $a$ is an increasing function of the input $I_0$: below $I_0=\theta$, no bump exists, and in any dimension, above $I_0=\theta$, a pulse exists. As the input amplitude increases, the pulse width increases. It is important to note that in this case that there exists a unique pulse for a given value of the input amplitude $I_0>\theta$. 

In both cases, we observe a quantitative and qualitative dependence of the zone of existence of pulses on the dimension of the space: while the size of the one-dimensional bump decreases almost linearly as $I_0$ is decreased, the size of the two dimensional bump has rather a quadratic behavior. This implies that there are different parameters corresponding to the existence of pulses: in the Mexican hat connectivity case, one needs larger inputs to obtain pulses (note that this ordering persists on a larger range of values tested, up to $I_0=50$), and at a given input amplitude value, and when both 1d and 2d models display the presence of a pulse, the small pulse is larger in the 2d model while the large pulse is smaller. In the inverse Mexican hat model, for small values of $I_0$ the two-dimensional pulse has a smaller half-width than the one-dimensional one, and for large values of $I_0$, it is the contrary. The question that arises is precisely what is the stability of these solutions. This is the topic of the following section. 

\subsection{Linear stability}
In order to characterize the linear stability of the pulses with respect to parameters, we need to characterize the spectrum of the linearized operator around these solutions. This is not a trivial task, because of the non-regularity of the function $f$. In order to handle such issues, Amari expressed the linearized operator in one dimension. This was then generalized by Bressloff and Kilpatrick~\cite{bressloff-kilpatrick:11} for two-dimensional models with synaptic depression and no delays. We extend these linearizations to the $d$-dimensional case with delays, and use a formalism based on the linearization in the sense of distribution instead of developments at leading orders to obtain these formulae. 

In order to characterize the linear stability of the bump $\mathcal{U}(\mx)$ with half-width $a$, we analyze the spectrum of the linearized operator around the bump. The method used by~\cite{bressloff-kilpatrick:11} consists in considering an arbitrary perturbation of small amplitude of the pulse $u(\mx,t) = \mathcal{U}(\mx)+\varepsilon {\varphi}(\mx,t)$, injecting this form of solution in the equation in order to find the evolution of the perturbation in the limit  $\varepsilon\to 0$:
\[\partial_t \varphi = \mathcal{L} \varphi\]
where
\[\mathcal{L}\varphi(\mx)=-\varphi(\mx) + \lim_{\varepsilon\to 0} \frac 1 {\varepsilon} \int_{\R^d} \J_d(\| \mx-\my\|) [f(\U(\my))-f(\U(\my)+\varepsilon \varphi(\my,t-\tau(\mx-\my)))]d\my.\]
A formal method allows to find this result in a more compact way using the differential of the Heaviside function $f$ in the sense of distributions, which is the Dirac delta function at $\theta$, $\delta_{\theta}(\cdot)=\delta(\cdot -\theta)$. This allows to write the linearized operator as:
\[\mathcal{L}\varphi(\mx) = -\varphi(\mx) + \int_{\R^d} \J_d(\|\mx-\my\|) \delta_{\theta}(\U(\my)) \varphi(\my,t-\tau( \mx-\my))) \; d\my. \]
The measure $\delta_{\theta}$ charges only the sphere of radius $a$ in $\R^d$ thanks to the rotational symmetry of $\U$ and by definition of the half-width of the pulse. In circular coordinates, denoting $\my= r \psi$ with $r\in \R_+$ and $\psi\in \S_{d-1}$ the unit sphere and using a change of variables, one readily obtains in the case of our particular bump (crossing once the threshold with spherical symmetry):
\[\mathcal{L}\varphi (\mx)= -\varphi (\mx)+ \frac{a^{d-1}}{\vert \U' (a)\vert}\int_{\S_{d-1}} \J_d( \|\mx -a\psi\|) \varphi(a \psi,t-\tau( \mx-a\psi))) \; {d_{d-1}\psi}\]
where $\vert \U' (a)\vert $ is the constant value taken by the gradient of $\U$ on the sphere of radius $a$ (see above).

In one dimension, the operator is therefore similar to the classical Amari linearized operation (with an additional delay term):
\[\mathcal{L}\varphi(x) = -\varphi(x) + \frac{{1}}{\vert \U'(a) \vert} \Big(\J_1(|x-a|) \varphi(a,t-\tau(x-a)) +\J_1(|x+a|) \varphi(-a,t-\tau(x+a)) \Big)\]
and in dimension 2, we recover the form of the linearized operator in~\cite{bressloff-kilpatrick:11}.

This linear operator includes a convolution term and a delay term. At this level of generality, it is not possible to determine a complete set of solutions to this equation. However, we may simplify the problem by considering complex perturbations in the class of functions $\varphi(\mx,t)=e^{\lambda t} p(\mx)$ for $\lambda \in \C$ and $p(\mx)$ a bounded continuous function on $\R^d$ decaying to zero at infinity (reflecting the properties of the weight distribution $\J_d$)\footnote{From the point of view of time evolution, the assumption is justified by noting that the action of the linear operator on the time variable is a convolution, which is therefore diagonalizable in the Fourier space. This is also the case of the space operator, which will motivate to specify, in a forthcoming section, to take an exponential shape for $p(\mx)$.}. Of course, the sign of the real part of $\lambda$ will give us the stability of the bump.

Let us now analyze the possible eigenvalues of the linear operator. Acting on the class of functions we considered, we rewrite the evolution in time of the perturbation and obtain a relationship, the \emph{dispersion relationship}, characterizing the eigenvalues $\lambda$:
\begin{equation}\label{eq:dispersionRelationship}
	(\lambda+1) p(\mx) = \frac{a^{d-1}}{\vert \U' (a)\vert} \int_{\S_{d-1}} \J_d( \| \mx -a\,\psi\|) e^{\lambda \tau( \mx-a\,\psi)} p(a\psi)\,d_{d-1}\psi.
\end{equation}
which rewrites, using the form of the delays considered:
\begin{equation*}
	(\lambda+1)e^{\lambda \tau_D} p(\mx) = \frac{a^{d-1}}{\vert \U' (a)\vert} \int_{\S_{d-1}} \J_d( \| \mx-a\,\psi\|) e^{\frac{\lambda}{c} \| \mx-a\,\psi\|} p(a\psi)\,d_{d-1}\psi.
\end{equation*}

The spectrum of the linear operator is composed of the complex values $\lambda$ such that the equation~\eqref{eq:dispersionRelationship} has solutions. This spectrum is comprised of the \emph{essential spectrum} corresponding to cases where the righthand side of the relationship vanishes (i.e. any perturbation vanishing on the sphere of radius $a$), and the \emph{point spectrum}.

The essential spectrum is therefore given by the solutions to the equation:
\[\lambda+1=0.\]
This corresponds to a stable mode of the bump: no destabilization may arise from perturbations vanishing on the sphere of radius $a$, and any such perturbation decays exponentially in time towards the bump solution. In other words, may the pulse be destabilized, the position of the threshold crossing will vary. 

The point spectrum is slightly more complex to handle, as it requires to take into account the infinite-dimensional operator on the righthand side of the dispersion relationship. It is precisely in this part of the spectrum that delays come into play and may produce qualitative effects such as destabilization of bumps. 

In order to illustrate the role of delays in the stability of the bump solutions, we now treat the case of a one-dimensional neural field analytically, before dealing with the problem in general dimension.

\subsubsection{One-dimensional neural fields}
The dispersion relationship greatly simplifies for one-dimensional neural fields allowing analytical treatment of the problem. In detail, in one space dimension, the dispersion relationship now reads:
\begin{equation}\label{eq:dispersionRelationship1d}
e^{ \lambda \tau_D} (\lambda + 1) p(x) = \frac{\J_1 (|x-a|)}{|\mathcal{U}' (a)|} p(a)e^{-\lambda \frac{|x-a|}{c}}+\frac{\J_1 (|x+a|)}{|\mathcal{U} ' (a)|} p(-a) e^{-\lambda \frac{|x+a|}{c}},
\end{equation}
and we search for complex solutions $\lambda$ of these equations. In the absence of delays, this is a linear equation in $\lambda$, and estimating this relationship at $x=\pm a$ yields a system of two equations which can be solved for $\lambda$ under the condition that $p(\pm a)\neq 0$. In the presence of delays, the equation is no more linear (it is a transcendental equation), but the same strategy allows obtaining analytical characterization of the eigenvalues. In detail, estimating the dispersion relationship at $x=a$ and $x=-a$ leads to the system:
\begin{equation*}
\left\{ 
\begin{array}{l}
(\J_1(0)-F(\lambda)) p(a) +  \J_1(2a)e^{-2 \frac{\lambda a}{c} }  p(-a)  = 0  \\ 
\J_1(2a)e^{-2 \frac{\lambda a}{c} } p(a) + (\J_1(0) - F(\lambda))  p(-a) = 0 \\
\end{array} 
\right.
\end{equation*}
where $F(\lambda) = e^{ \lambda \tau_D}( \lambda + 1 ) | \mathcal{U} ' (a)|$.  
This system is linear in $p(\pm a)$. It has has non trivial solutions for $p( \pm a)$ if and only if its determinant vanishes:
\begin{equation}\label{eq:CondBif}
(\J_1(0)-F(\lambda))^2 = \J_1(2a)^2e^{-4 \frac{\lambda a}{c} }
\end{equation} 
which leads to the relationship: 
\begin{equation*}
(\lambda_{\pm}+1)e^{ \lambda_{\pm} \tau_D} = \frac{\J_1(0) \pm \J_1(2a)e^{-2 \frac{\lambda_{\pm} a}{c} }}{| \mathcal{U} ' (a)|}.
\end{equation*}
Once this relationship is solved, we obtain the eigenfunctions through the dispersion relationship~\eqref{eq:dispersionRelationship1d}:
\[ p_{\pm}(x) = \frac{p(a)}{F(\lambda_{\pm})} \left(\J_1(|x-a|)e^{-\lambda_{\pm} \frac{|x-a|}{c}} +  \J_1 (|x+a|) e^{-\lambda_{\pm} \frac{|x+a|}{c}}\frac{F(\lambda_{\pm})-\J_1(0)}{\J_1(2a)}e^{\frac{2\lambda_{\pm} a}{c}}\right)\]
and using the relationship~\eqref{eq:CondBif}, we know that the coefficient multiplying the second term within the parenthesis has a square equal to one, hence we have two possible solutions:
\begin{equation}\label{eq:DestabilizationModes}
	p_{\pm}(x) = \frac{p(a)}{F(\lambda_{\pm})} \left(\J_1(|x-a|)e^{-\lambda_{\pm} \frac{|x-a|}{c}} \pm  \J_1 (|x+a|) e^{-\lambda_{\pm} \frac{|x+a|}{c}}\right).
\end{equation}

\begin{figure}
	\centering
		\subfigure[Symmetric Perturbation]{\includegraphics[width=.3\textwidth]{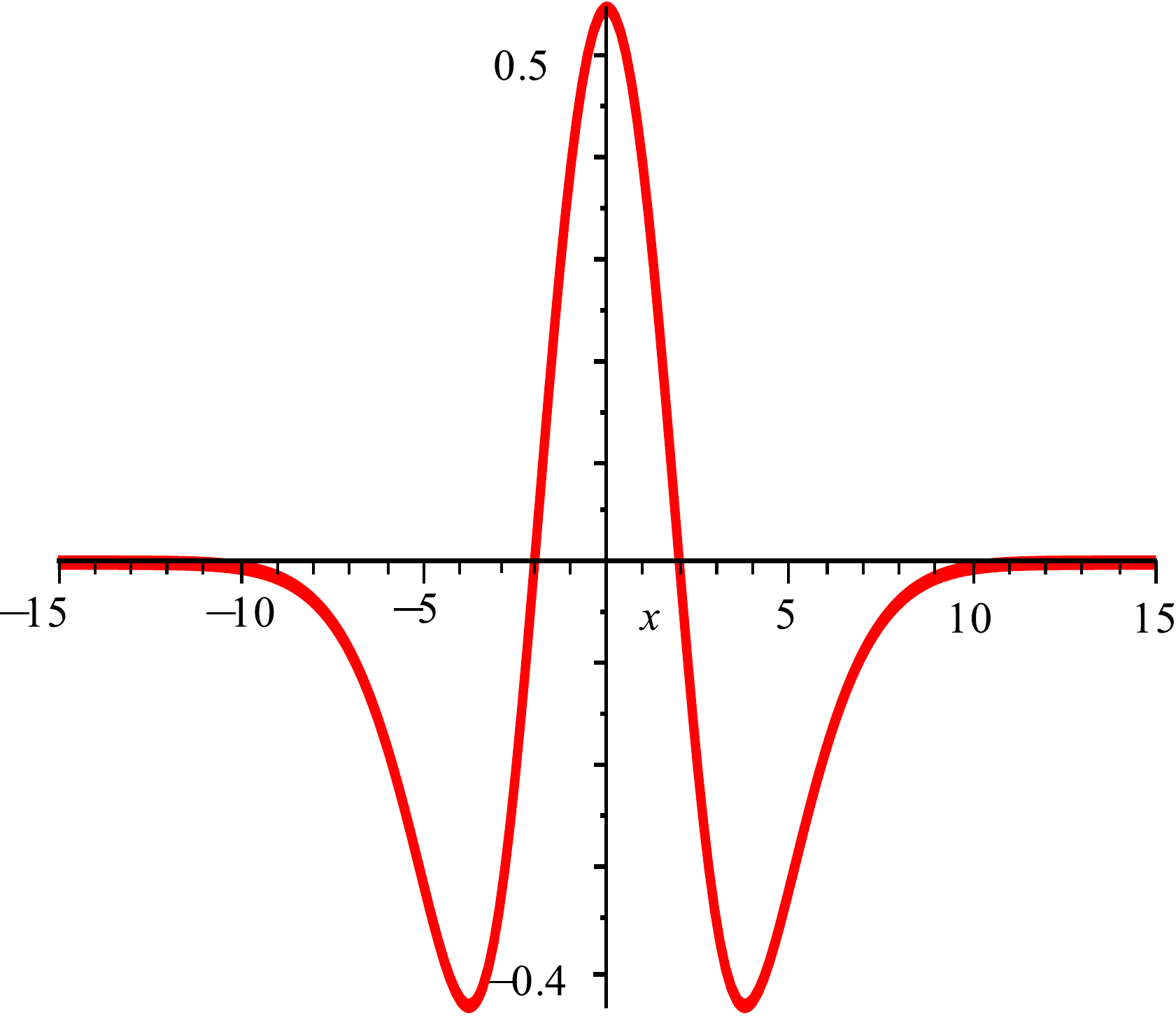}}\qquad\quad\quad
		\subfigure[Asymmetric Perturbation]{\includegraphics[width=.3\textwidth]{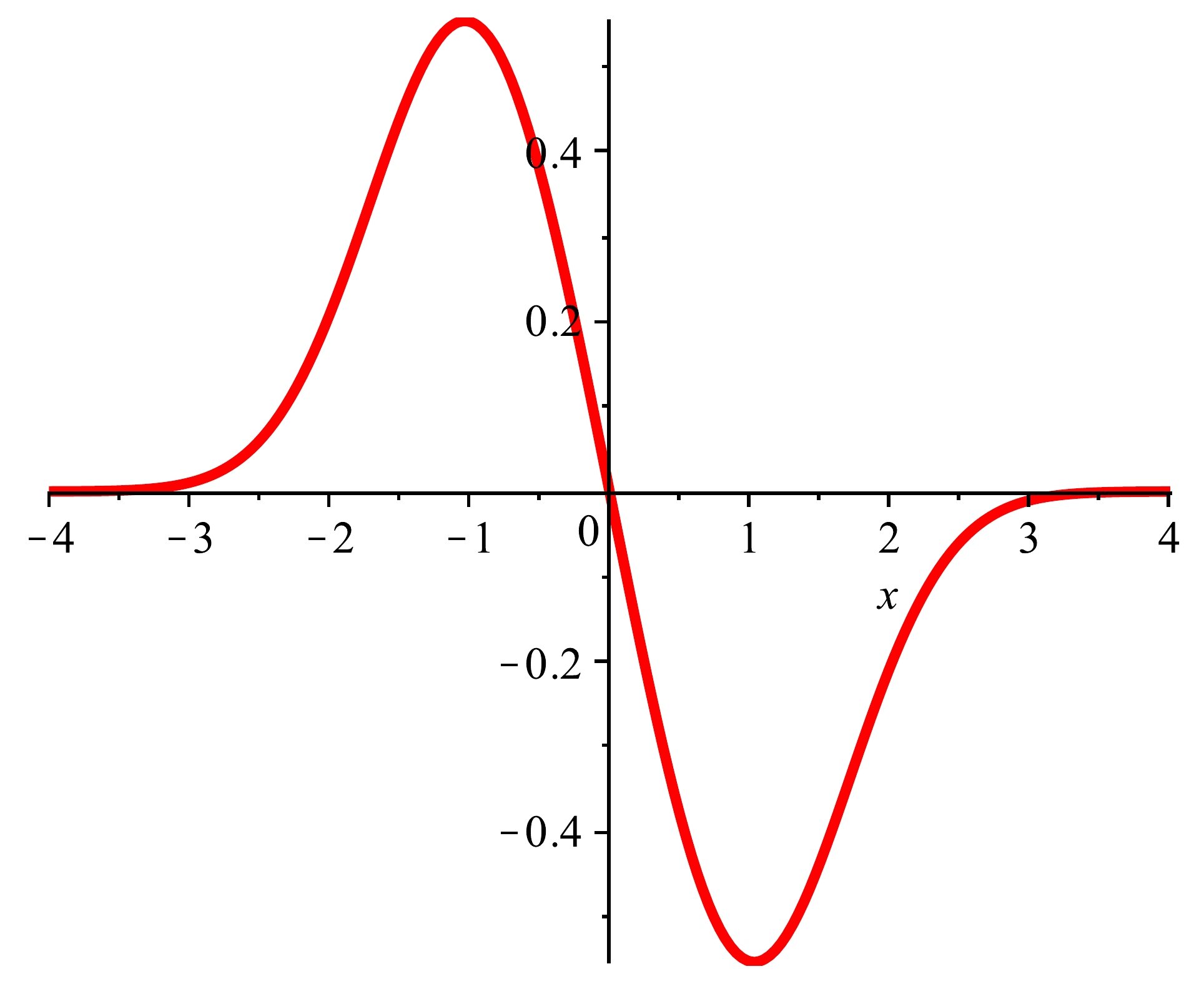}}
	\caption{Form of the eigenfunctions related to the bump destabilization for the inverse Mexican hat (Model I, see Fig~\ref{fig:Bump1DSynapticDelayModelI} for the parameters): (a) symmetric perturbation $p_+$ corresponding to breathing waves and (b) asymmetric perturbations $p_-$ corresponding to sloshing waves.}
	\label{fig:Perturbations}
\end{figure}

We note that these are complex-valued functions. They have to be seen as a two-dimensional eigenspace spanning the vector space corresponding to the two functions $\Re [e^{\lambda_{\pm} t}p_{\pm}(x)]$ and $\Im [e^{\lambda_{\pm} t}p_{\pm}(x)]$ (and correspond to the two complex conjugate eigenvalues $\lambda_{\pm}$ and $\lambda_{\pm}^*$). Interestingly, the two eigenmodes found clearly correspond to a symmetric and asymmetric perturbation of the bump, and an instability along one of these modes allows to infer the type of pattern observed. These two eigenfunctions are plotted in Fig.~\ref{fig:Perturbations} for a specific choice of parameter. 

We are interested in changes in the stability of the bump solution. These changes occur when one can find a characteristic root with zero real part. These are of two types: either fold transitions when $\lambda=0$, or Hopf transitions when the characteristic roots are purely imaginary. Fold transitions are independent of the delays. In our model, they therefore arise exactly for the same parameters as in the non-delayed equation, and were largely studied in the literature (see e.g.~\cite{amari:77,F}). Analysis shows that for delays equal to zero, the stability of the pulses can be related to a geometric property of the curves $I_0(a)$ corresponding to the existence of pulses (see Fig.~\ref{fig:I0}): one can show that pulses are stable whenever $I_0'(a)>0$ and unstable otherwise, justifying the stability depicted in the diagram. 

The question we may now ask from this setting is whether delays induce changes in the stability of these bumps. From what precedes, it is clear that such destabilization necessarily occur through Hopf bifurcations. In order to characterize such bifurcations, we solve equation~\eqref{eq:dispersionRelationship} for purely imaginary solutions $\lambda_{\pm} = i \omega_{\pm}$. Necessarily, such solutions exist if and only if the modulus of the righthand side is greater than zero, and in that case, equating modulus and argument in the previous relationship yields two conditions to be satisfied in the parameter space so that Hopf bifurcations occur.

Let us first address the case of constant delays (we start by neglecting the propagation delay term, i.e. considering $c \to \infty$). Hopf bifurcations arise if and only if
\[\left|\frac{\J_1(0) \pm \J_1(2a)}{| \mathcal{U} ' (a)|} \right|> 1,\]
and for parameters satisfying the simple relationship:
\begin{equation}\label{eq:Bifurcation1d}
\left\{ 
\begin{array}{l}
\omega_{\pm} = \sqrt{  \left (  \frac{\J_1(0) \pm \J_1(2a)}{| \mathcal{U} ' (a)|} \right  )^2  -1}  \\ 
{\tau_{D}}_{\pm} \equiv \frac{1}{\omega_{ \pm}} \left ( \arg \left ( \frac{\J_1(0) \pm \J_1(2a)}{| \mathcal{U} ' (a)|} \right )-\text{arccos}\left(\frac{| \mathcal{U} ' (a)|}{|\J_1(0) \pm \J_1(2a)|} \right)\right ) \pmod{\frac{2\pi}{\omega_\pm}}.
\end{array} 
\right.
\end{equation}
The above notation for the delays ${\tau_{D}}_\pm$ has to be understood has follows. If $\arg \left ( \frac{\J_1(0) \pm \J_1(2a)}{| \mathcal{U} ' (a)|} \right )\geq\text{arccos}\left(\frac{| \mathcal{U} ' (a)|}{|\J_1(0) \pm \J_1(2a)|} \right)$, then we set the critical delay to
\bqs
{\tau_{D}}_{\pm}= \frac{1}{\omega_{ \pm}} \left ( \arg \left ( \frac{\J_1(0) \pm \J_1(2a)}{| \mathcal{U} ' (a)|} \right )-\text{arccos}\left(\frac{| \mathcal{U} ' (a)|}{|\J_1(0) \pm \J_1(2a)|} \right)\right )\geq 0,
\eqs
otherwise, we set it to
\bqs
{\tau_{D}}_{\pm}= \frac{1}{\omega_{ \pm}} \left (2\pi+ \arg \left ( \frac{\J_1(0) \pm \J_1(2a)}{| \mathcal{U} ' (a)|} \right )-\text{arccos}\left(\frac{| \mathcal{U} ' (a)|}{|\J_1(0) \pm \J_1(2a)|} \right)\right )\geq 0.
\eqs

Formula \eqref{eq:Bifurcation1d} provides a closed-form relationship, in the parameter space, where Hopf bifurcations occur. This relationship depends on half-width $a$ of the bump, which is a function of other parameters given by equation~\eqref{eq:BumpExistence}, and as we saw, this can be written in a simple form when solved for $I_0$ (considering $a$ as a parameter). This allows to have a closed-form relationship between parameters corresponding to the Hopf bifurcation curve, along the symmetric or asymmetric mode. 

\begin{figure}
	\centering
		\includegraphics[width=.7\textwidth]{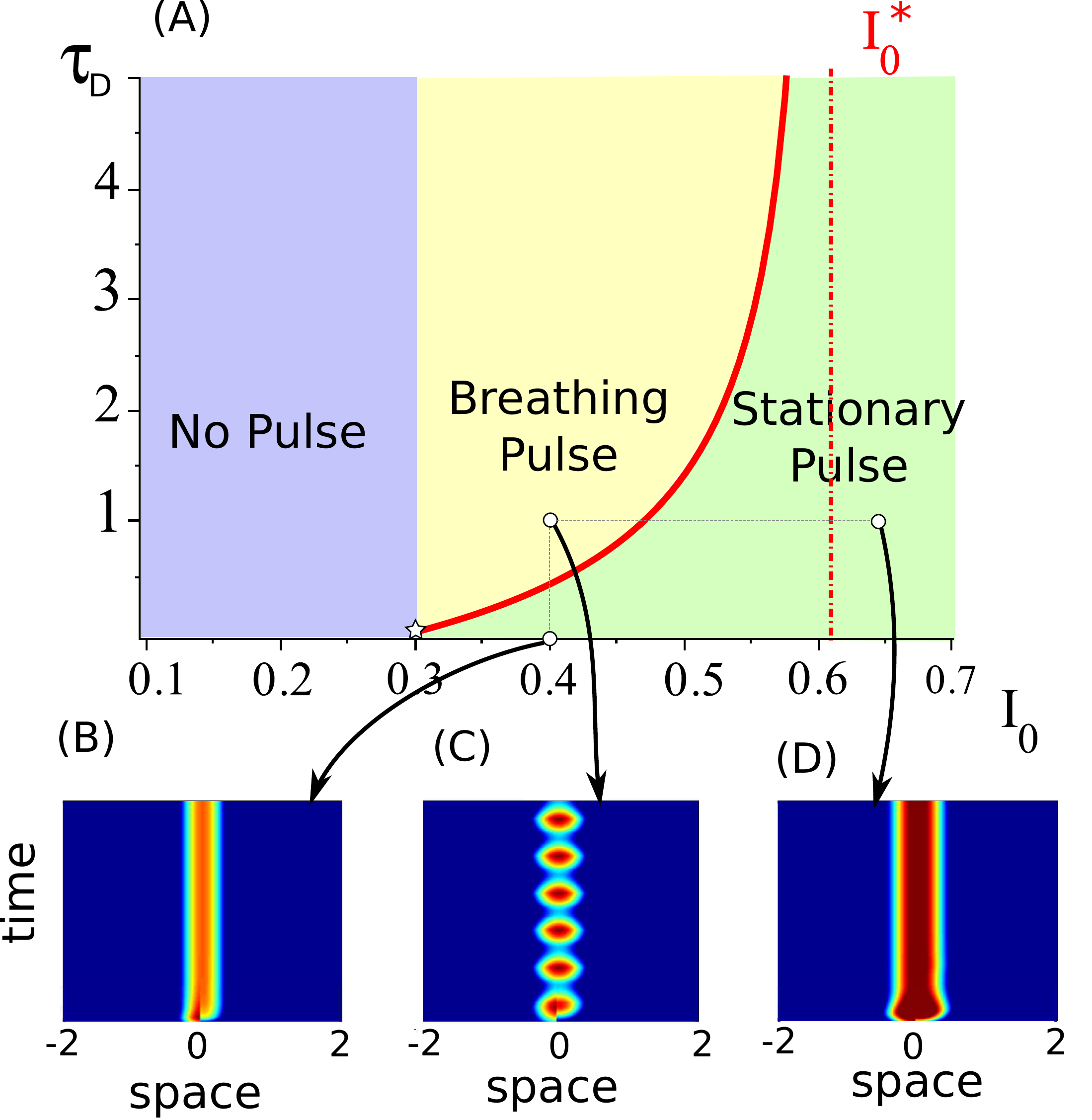}
	\caption{Delay-induce breathing pulse in one dimensional neural field with inverse Mexican hat connectivity (Model I): $\omega_e=1.3$, $\omega_i=1.1$, $\sigma_e=4$, $\sigma_i=2$, $\theta=0.3$, $\sigma=1.5$. (A): codimension 2 bifurcation diagram in the plane $(\tau_D,I_0)$. Blue region: no pulse ($I_0\leq \theta$), Yellow: breathing pulses and Green: Stationary pulses. Red curve corresponds to the Hopf bifurcation curve~\eqref{eq:Bifurcation1d}, and the red dashed-dotted line to the value of $I_0^*$. (B), (C) and (D) are space-time representations of the solution for different choices of parameters: (B) $I_0=0.4$, $\tau_D=0$: stationary pulse, (C): $I_0=0.4$, $\tau_D=1$: breathing pulse, (D): $I_0=0.65$, $\tau_D=1$: stationary pulse. Colormap with fixed boundaries $\theta \pm 20\%$}
	\label{fig:Bump1DSynapticDelayModelI}
\end{figure}

\begin{figure}
	\centering
		\includegraphics[width=\textwidth]{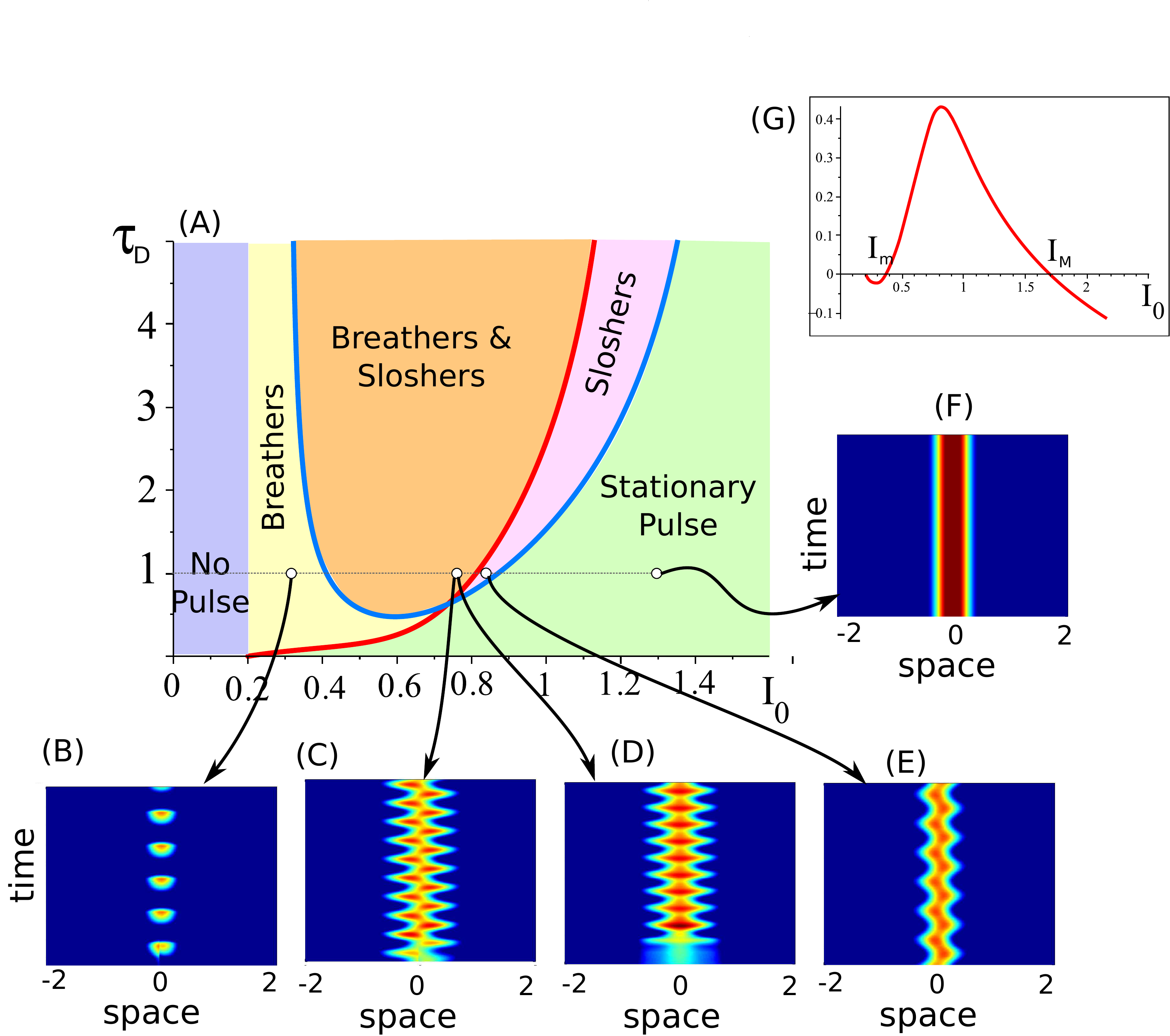}
	\caption{Delay-induce breathers and sloshers in one dimensional neural field with inverse Mexican hat connectivity (Model II): $\omega_e=1$, $\omega_i=1.5$, $\sigma_e=1.5$, $\sigma_i=1$, $\theta=0.2$, $\sigma=1$. (A): codimension 2 bifurcation diagram in the plane $(\tau_D,I_0)$. Blue region: no pulse ($I_0\leq \theta$), Yellow: symmetric mode unstable, asymmetric mode stable, pink: symmetric mode stable, asymmetric mode unstable, orange: both symmetric and asymmetric modes unstable, and Green: Stationary pulses (no instability). Red (resp. blue) curve corresponds to the Hopf bifurcation along the symmetric (resp. asymmetric) direction given by equations~\eqref{eq:Bifurcation1d}. (B-F) are space-time representations of the solution for $\tau_D=1$ and different choices of input amplitude $I_0$: (B) $I_0=0.3$: breathing pulse, (C-D): $I_0=0.75$, with a symmetric initial condition (C) and asymmetric initial condition (D), (E): $I_0=0.82$, symmetric initial condition: sloshing wave, and (F) $I_0=1.2$: stationary pulse. (G) represents the map $I_0\mapsto 2(\J_1(2a)-\J_1(0))-I'(a)$ with $a$ the half-width of the pulse associated with $I_0$. In the region where this map is positive, Hopf bifurcations along the asymmetric mode arise, defining the quantities $I_m$ and $I_M$.}
	\label{fig:Bump1DSynapticDelayModelII}
\end{figure}

More precisely, we have:
\renewcommand{\theenumi}{(\roman{enumi})}
\begin{enumerate}
	\item If $\left|\frac{\J_1(0) + \J_1(2a)}{| \mathcal{U} ' (a)|}\right| > 1$, a Hopf bifurcation arises when $\tau_D$ reaches the value given in equation~\eqref{eq:Bifurcation1d}, and the destabilization arises along a symmetric mode $p_+(x)$ given by equation~\eqref{eq:DestabilizationModes}. This gives rise to localized periodic oscillations symmetric around the stationary pulse, or \emph{breather} (see e.g. Fig.~\ref{fig:Bump1DSynapticDelayModelI}(C) or Fig.~\ref{fig:Bump1DSynapticDelayModelII}(B,D)).
	\item If $\left|\frac{\J_1(0) - \J_1(2a)}{| \mathcal{U} ' (a)|} \right|> 1$, a Hopf bifurcation arises when $\tau_D$ reaches the value given in equation~\eqref{eq:Bifurcation1d}, and the bump looses stability along the antisymmetric mode $p_-(x)$. This gives rise to an asymmetric periodic solution around the stationary pulse, a \emph{slosher} (see Fig.~\ref{fig:Bump1DSynapticDelayModelII}(C,E)).
\end{enumerate}

Which mode is unstable, and therefore the type of periodic solution obtained passed the Hopf bifurcation depends on the choice of parameters, and particularly on the type of connectivity chosen.

For Mexican hat connectivity, it is folklore that no delay-induced destabilization is likely to arise. This has been suggested in  a number of studies~\cite{roxin-montbrio:11,veltz-faugeras:11,veltz2013}, and is evocative of the analogous fact, valid for ordinary differential equations in one dimension, that only inhibitory delayed feedbacks loops can destabilize a fixed point. In the present setting, the same phenomenon occurs: stable pulses do not undergo Hopf bifurcations. The unstable pulse does undergo a Hopf bifurcation along the sum mode (symmetric perturbation), but the dynamics of the system remain qualitatively the same: the pulse with largest half-width is globally stable, while the pulse with smallest half-width remains unstable whatever the delays.

The inverse Mexican hat model is much more sensitive to the presence of delays. We now focus on two examples. The first example is chosen so that the total amplitude of the connectivity function is small enough, preventing possible destabilization along the asymmetric mode. The second case considers a large inverse Mexican hat connectivity, and the asymmetric mode will be destabilized. 

In Fig.~\ref{fig:Bump1DSynapticDelayModelI}, we consider a first choice of inverse Mexican hat (Model (I)) whose positive and negative peaks are small enough, so that no instability arise from the asymmetric mode (no $\omega_-$ exist solving equations~\eqref{eq:Bifurcation1d}). In detail, parameters are chosen so that $2(\J_1(2a)-\J_1(0))<I'(a)$ for all $a$. In that case, any perturbations along the asymmetric mode cancel away, while the symmetric mode $p_+$ gets destabilized for sufficiently large delays in a specific region of input amplitude $I_0$. This region is bounded on the left by the value $I_0=\theta$ corresponding to the emergence of pulses (condition~\eqref{eq:BumpExistence}), and on the right by a value $I_0^*$ corresponding to the condition:
\[I_0^*=-\J_1(2a)\frac{\sigma^2}{a} e^{\left(\frac{a}{\sigma}\right)^2},\]
for which $\J_1(0)-\J_1(2a) + \I'(a)=\J_1(0)+\J_1(2a)$, hence $\omega_+\to 0$ and $\tau \to \infty$. Note that this value necessarily exists since as $I_0$ increases, $a$ increases and $I'(a)$ decreases towards zero, while $\J_1(2a)$ is positive for $a$ sufficiently large. At the emergence of pulses (for $I_0=\theta$), we are in a situation such that $a=0$ and $\U'(a)=0$ (tangent case), in which case $\omega\to \infty$ and $\tau_D\to 0$. For any value of the input amplitude $I_0\in [\theta, I_0^*]$, there exists a value of $\tau_D$ for which the system undergoes a Hopf bifurcation and the symmetric mode $p_+$ is destabilized. This gives rise to oscillations around the stationary pulse along the symmetric mode, i.e.,  breathers. Interestingly, sufficiently large input ($I_0>I_0^*$) yield unconditionally stable stationary pulses, that cannot be destabilized by the presence of delays. 

For inverse Mexican hat connectivities with larger amplitude, one can find real values $\omega_-$ solutions to equation~\eqref{eq:Bifurcation1d}, and therefore asymmetric instabilities (sloshing pulses, or sloshers). This is the case of Model (II) displayed in Fig~\ref{fig:Bump1DSynapticDelayModelII}. The breathing instability persists. In that model, we therefore obtain two curves of Hopf bifurcations, one corresponding to the destabilization along the symmetric mode $p_+$ (red curve) and the other to the asymmetric mode $p_-$ (blue curve). The symmetric mode destabilization curve has the same properties as those described for model (I): it emerges from $\tau_D=0$ at a value $I_0=\theta$ corresponding to the existence of pulses condition, and diverges for a finite value of $\tau_D$ corresponding to $I_0^*$. The Hopf bifurcation corresponding to the asymmetric mode does not have the same shape: indeed, in order to obtain a solution $\omega_-$ to equations~\eqref{eq:Bifurcation1d}, one needs to exploit the fact that the connectivity does not have a constant sign, and therefore the Hopf bifurcation arises for values of $I_0$ corresponding to sufficiently large pulse half-width. These exist when 
\[2(\J_1(2a)-\J_1(0)) \geq \I'(a)\]
which defines an interval of values $[a_m,a_M]$, or of input amplitudes $I_0\in[I_m,I_M]$ for which the condition is satisfied (see Fig~\ref{fig:Bump1DSynapticDelayModelII}(G)). At the boundaries of this interval, as discussed in Model (I), the value of $\omega_-$ tends to zero, and therefore the delay $\tau_D$ corresponding to the Hopf bifurcation diverges to infinity. Let us eventually note that parameter domains for which the destabilization along the symmetric and asymmetric modes arise do not completely overlap: it is we observe that $I_m>\theta$ and $I_M>I_0^*$. These curves therefore partition the parameter space $(I_0,\tau_D)$ into 5 regions: no pulse, stationary pulse, breathers (instability of the symmetric mode only), sloshers (instability of the asymmetric mode only), and a mixed breathers-slosher regions where both symmetric and asymmetric modes are unstable. 

Let us now investigate the role of propagation delays. We now no more assume that the speed of propagation of information is instantaneous. We denote 
\[\Phi(\tau_D,\omega,a) = \frac{\J_1(0)-\vert \U'(a)\vert e^{\mathbf{i}\omega \tau_D}(\mathbf{i}\omega+1)}{\J_1(2a)^2}.\]
It is easy to see, equating modulus and argument in the dispersion relationship~\eqref{eq:CondBif}, that Hopf bifurcations arise when parameters satisfy the equations: 
\begin{equation}\label{eq:BifPropagDelays}
\begin{cases}
	|\Phi(\tau_D,\omega,a)|=1 \\
	c = - \frac{4 a \omega}{\arg (\Phi(\tau_D,\omega,a)^2) + 2 \pi m}
\end{cases}
\end{equation}
for some $m \in \mathbb{Z}$. These formulae are much less transparent than in the case with constant delays. However, these are relatively simple when investigating the respective roles of constant and propagation delays. In that case, for fixed connectivity, threshold and input, the width of the stationary pulse can be readily computed using equation~\eqref{eq:BumpExistence}, and it is easy to write the Hopf bifurcation curve given by equation~\eqref{eq:BifPropagDelays} as a parametric curve (parameterized by $\omega$). 
\begin{figure}[htbp]
	\centering
		\includegraphics[width=0.8\textwidth]{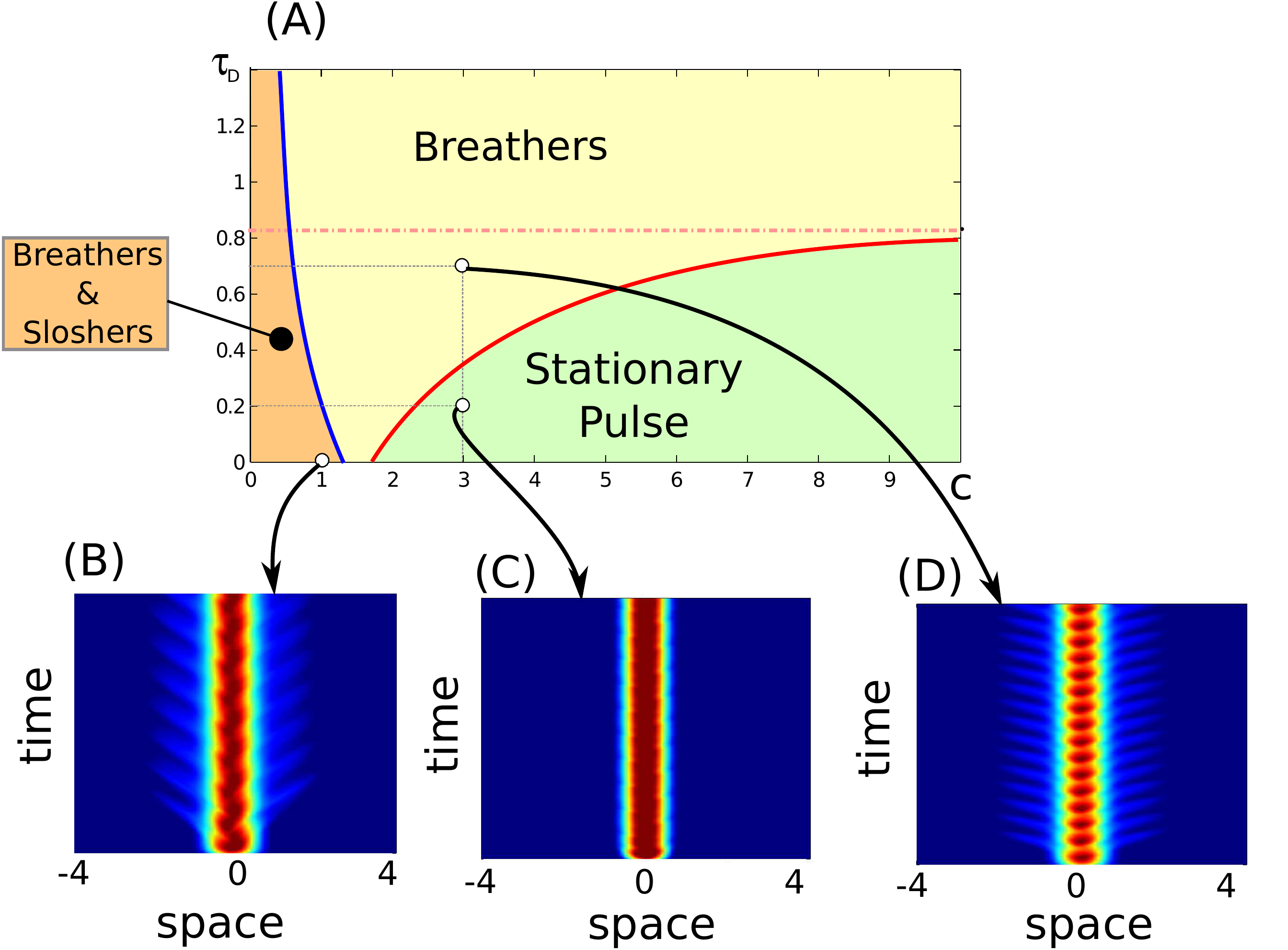}
	\caption{Stability of the pulse with propagation delays for Model (I) with $I_0=0.4$. (A) Bifurcation diagram in the $(c,\tau_D)$ plane with the first Hopf bifurcation (red line) that asymptotes the value of the delay corresponding to the purely synaptic case ($\tau_D\approx 0.815$, pink dot-dashed line). (B) $c=1$, $\tau_D=0$: sloshing pulse, (C): $c=3$, $\tau_D=0.2$: stationary pulse, (D) $c=3$, $\tau_D=0.7$: breathing pulse. }
	\label{fig:propagationDelay}
\end{figure}
This method allows to draw the Hopf bifurcation curve obtained in Figure~\ref{fig:propagationDelay}. The effective role of propagation delays on the stability of the pulse along the symmetric mode is relatively intuitive: indeed, by increasing the effective averaged delay, the finite-speed propagation delay shifts the bifurcation point to smaller values of the constant delay $\tau_D$. This explains the increasing shape of the bifurcation curve for Model (I) analyzed in Fig.~\ref{fig:Bump1DSynapticDelayModelI} with an asymptote at the value predicted by the analysis of the case with no propagation. Interestingly, we observe that below a certain speed, the stationary pulse is unstable even with no constant delay. More surprising is the emergence of a destabilization of the pulse along the asymmetric model at low propagation speed (while the asymmetric mode was unconditionally stable in the limit $c\to\infty$), leading to the possibility of sloshing waves induced by the shape of the propagation delays. Eventually, we note that the shape of the breathing and sloshing pulses reflects the presence of propagation delays. Indeed, we observe (see Fig.~\ref{fig:propagationDelay} (B) and (D)) that the oscillations of the shape of the pulse induce evanescent waves traveling in a small region around the pulse with a specific speed that depends on the propagation speed $c$, as visible from the more horizontal direction of these waves in the case $c=3$ (D) compared to $c=1$ (B).

\subsubsection{Pulses in general dimensions}
The one-dimensional neural field case is informative in that all formulae can be obtained in closed form, allowing to understand the role of parameters in the stability of the pulses, as well as possible instabilities along symmetric modes (breathers) or asymmetric modes (sloshers). In the one-dimensional context indeed, asymmetric perturbations are relatively simple: the instability arises along a single mode, $p_-$, which is an odd function. In higher dimensions, due to the rotational symmetry and geometric constraints, there exists an infinite number of asymmetric modes around the pulse. For radially symmetric input, pulses correspond to solutions $\U(\mx)$ with spherical symmetry, i.e. only depending on $r=\Vert \mx \Vert$, which we denote, with a slight abuse of notations $\mathcal{U}(r)$, and for which there exists $a\geq 0$ such that:
\begin{equation}
\begin{cases}
	\U(r) > \theta, & r<a \\ 
	\U(r) < \theta, & r > a\\
	\mathcal{U}(a) = \theta.
\end{cases}
\end{equation}
We recall that the pulse is known explicitly as a function of $\mathcal{M}$ and $I$:
\begin{equation}
\mathcal{U}(r)= \mathcal{M}(r,a) + I(r), 
\end{equation}
and the stability of the bump is given by the eigenvalues of the linearized operator, solutions of the general dispersion relationship~\eqref{eq:dispersionRelationship}. 

In order to characterize the stability of the pulse, we therefore need to solve this equation. The eigenvectors $p(\mx)$ provide the modes along which the pulse is deformed if such an instability is found. Let us start by considering radially symmetric functions $p(\mx)=p(r)$. In this case, he dispersion relationship greatly simplifies when evaluated on the sphere of radius $a$. One obtains the relationship:
\[	(\lambda+1)e^{\lambda \tau_D} = \frac{a^{d-1}}{\vert  \U'(a) \vert} \int_{\S_{d-1}} \J_d(a \| \mathbf{e}_1-\psi\|) e^{\frac{\lambda}{c} a\| \mathbf{e}_1-\psi\|}\,d_{d-1}\psi.
\]
where $\mathbf{e}_1$ is the first vector of the canonical basis (without loss of generality). In the absence of propagation delays ($c\to\infty$), we therefore obtain a relatively simple equation:
\[	(\lambda+1)e^{\lambda \tau_D} = \frac{a^{d-1}}{\vert  \U'(a) \vert} \int_{\S_{d-1}} \J_d(a \| \mathbf{e}_1-\psi\|)\,d_{d-1}\psi=:\Phi_d(a)
\]
similar to those obtained in one dimension. In particular, the righthand side of the equation is a simple function $\Phi_d(a)$ of the parameters of the model. The same algebra as done in one dimension yields the equations for purely imaginary solutions $\lambda = \mathbf{i}\omega$:
\begin{equation}\label{eq:Bifurcation_nd}
\begin{cases}
	\omega = \sqrt{\Phi_d(a)^2  -1}  \\ 
	\tau_{D} \equiv \frac{1}{\omega} \left (\gamma(a) -\text{arccos}\left(\frac 1 {\vert \Phi_d(a)\vert}\right) \right ) \pmod{\frac{2\pi}{\omega}}.
\end{cases}
\end{equation}
with $\gamma(a)=0$ if $\Phi_d(a)\geq 0$ and $\pi$ otherwise.

\begin{figure}[h]
	\centering
		\includegraphics[width=.8\textwidth]{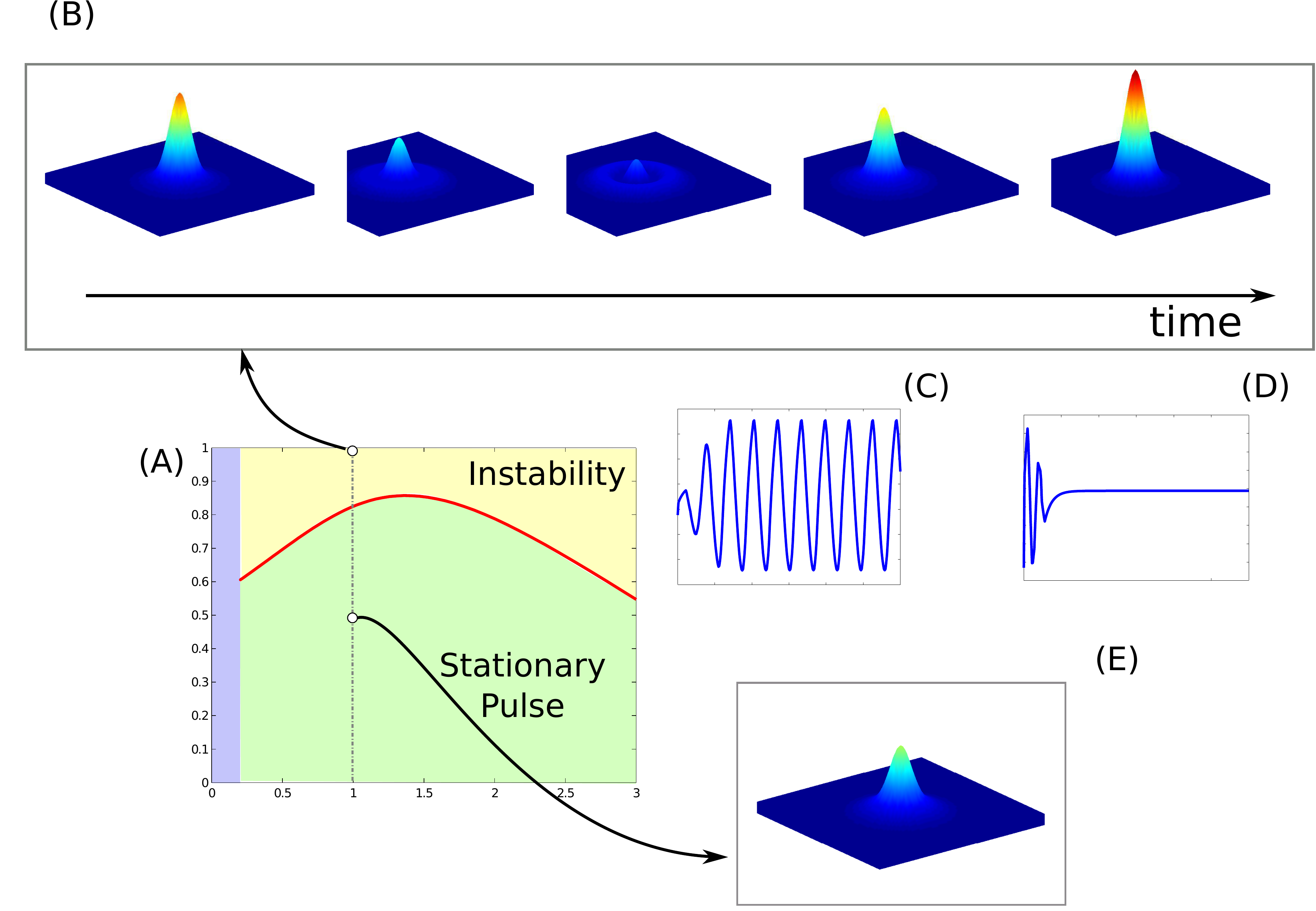}
	\caption{Breathers in two dimensions. Parameters: $w_e=2$, $w_i=2.5$, $\sigma_e=1$, $\sigma_i=0.5$, $\sigma=0.5$ and $\theta=0.3$. (A) represents the Hopf bifurcation curve along the symmetric direction (which was the only instability found in this system). The evolution in the instability region ($I_0=1$ and $\tau_D=1$) is provided in panel (B): five snapshots over one period of the breather were taken, and the time evolution of the amplitude at the center is plotted in panel (C). Analogous plots are provided in the region where no instability was found: the pulse does not evolve in time, as visible from the five snapshots and evolution of the central amplitude (D-E). Movies of the evolution are provided in supplementary material. }
	\label{fig:Bump2d}
\end{figure}

Asymmetric modes are more complex to characterize. However, due to the particular form of the linearized operator, we observe that the only values of the perturbation that are active are those on the $\R^d$-sphere of radius $a$. Thanks to the form of the equation, it is easy to see that the family of exponentials diagonalize the linear operator. Indeed, for $p(x)=p(0)e^{\mathbf{i} \mathbf{l}\cdot \psi}$ with $\mathbf{l}\in \R^d$, we obtain, evaluating the dispersion relationship at $a\mx$ on the sphere of radius $a$ ($\mx$ is on the unit sphere):
\begin{align*}
	e^{\lambda \tau_D}(\lambda+1)e^{\mathbf{i} \mathbf{l}\cdot a \mx} &= \frac{a^{d-1}}{\vert \U'(a)\vert} \int_{\S_{d-1}}\J_d(a \| \mx-\psi\|)\,e^{\mathbf{i} a \mathbf{l}\cdot \psi}d_{d-1}\psi\\
	&=\frac{a^{d-1}}{\vert \U'(a)\vert} \int_{\mx+\S_{d-1}}\J_d(a \|\mathbf{z}\|)\,e^{-\mathbf{i} a \mathbf{l}\cdot \mathbf{z}}d_{d-1}\mathbf{z}\; e^{\mathbf{i} \mathbf{l}\cdot a \mx}
\end{align*}
yielding the dispersion relationship:
\[e^{\lambda \tau_D}(\lambda+1) = \Phi_d(a,\mathbf{l})\]
with 
\[\Phi_d(a,\mathbf{l}) = \frac{a^{d-1}}{\vert \U'(a)\vert} \int_{\mathbf{e}_1+\S_{d-1}}\J_d(a \|\mathbf{z}\|)\,e^{-\mathbf{i} a \mathbf{l}\cdot \mathbf{z}}d_{d-1}\mathbf{z}.\]
Here, we used the spherical symmetry to state that the integrals on the righthand side are all equal whatever $\mx$, and denoted $\mathbf{e}_1$ the first vector of the canonical basis, which belongs to $\S_{d-1}$. This dispersion relationship is solved exactly as done in equation~\eqref{eq:Bifurcation_nd}, providing the possible Hopf bifurcations along asymmetric modes (non constant on the sphere of radius $a$).

These formulae generalize the approach developed in~\cite{folias-bressloff:04} in two-dimensional neural fields. In this case, the sphere is the unit circle, and the dispersion relationship reads: 
\begin{equation*}
e^{ \lambda \tau_D}(\lambda + 1) p (\textbf{r}) =  \frac{a}{|\mathcal{U}'(a)|} \int\limits_{0}^{2 \pi} \J_2(\|\textbf{r}-a \, e^{\mathbf{i}\theta}\|) p(a e^{\mathbf{i}\theta}) d \theta.
\end{equation*}
Estimating this relationship at $\textbf{r}=\textbf{a}$, one obtains the condition:
\begin{equation}
e^{ \lambda \tau_D}(\lambda + 1) p (\theta) = \frac{a}{|\mathcal{U}'(a)|} \int\limits_{0}^{ \pi} \J_2 \left (2 a |\sin \left (\theta - \theta' \right ) | \right ) p(2 \theta') d \theta '
\end{equation}
Now we choose angular perturbation of the form $p(\theta) = e^{i n \theta}$ to get : 
\begin{equation}\label{eq:HigherOrder}
e^{ \lambda \tau_D}(\lambda + 1) = \frac{a}{|\mathcal{U}'(a)|} \int\limits_{0}^{ \pi} \J_2(2 a |\sin( \theta)|) \cos(2 n \theta) d \theta =: J_n(a).
\end{equation}

\begin{figure}[ht]
	\centering
		\includegraphics[width=\textwidth]{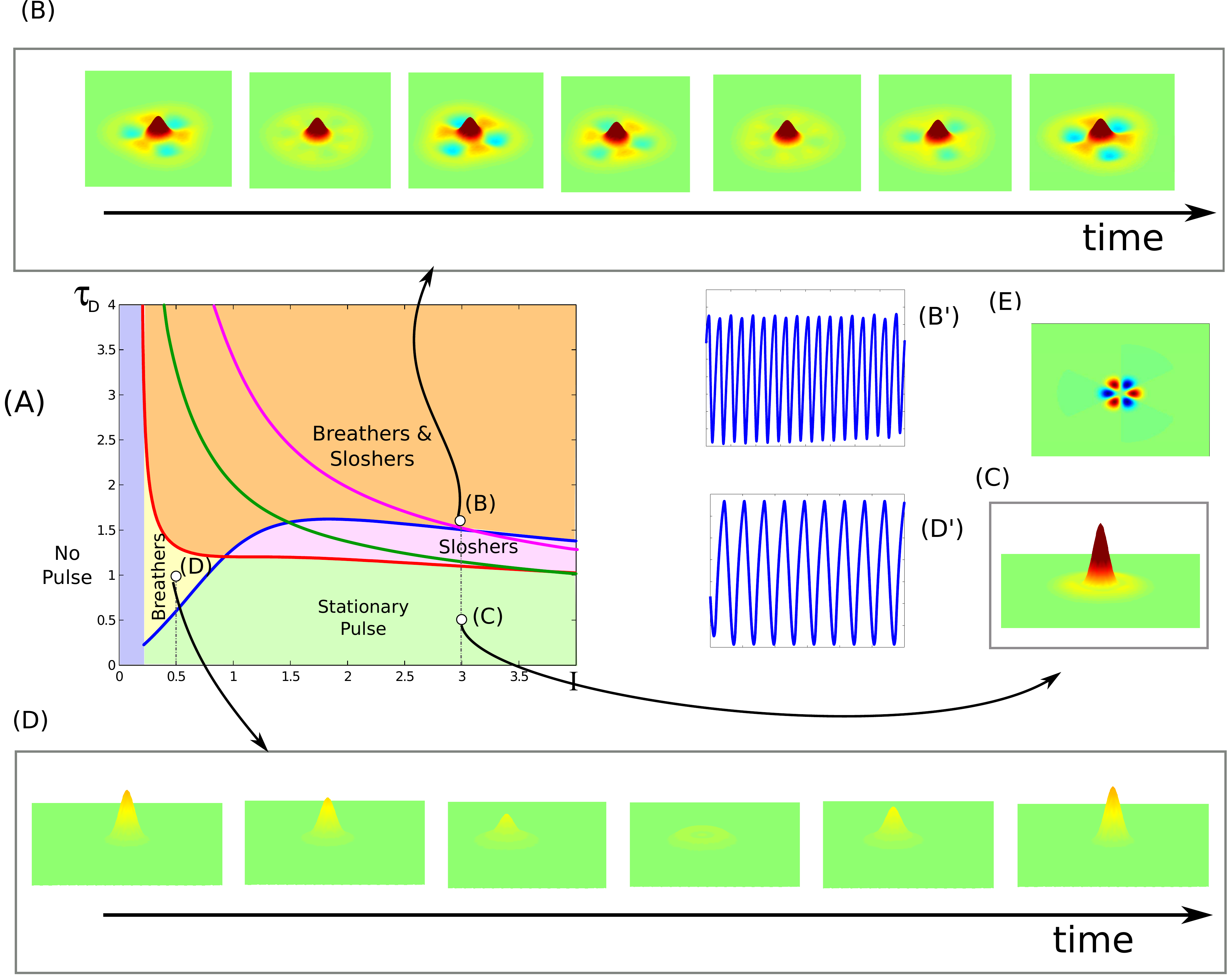}
	\caption{Breathers and sloshers in two dimensions. Parameters: $w_e=2.4$, $w_i=2$, $\sigma_e=1$, $\sigma_i=0.5$, $\sigma=0.5$ and $\theta=0.2$. (A) represents the Hopf bifurcation curve along the symmetric direction (blue curve) together with the first three instability modes ($n=1$, blue, $n=2$, green, $n=3$, magenta). These lines, similarly to the case of Fig~\ref{fig:Bump1DSynapticDelayModelII}, decompose the parameter space into 5 regions (same color code as in Fig.~\ref{fig:Bump1DSynapticDelayModelII}). Initial conditions were perturbations of the stationary pulse (found for $\tau_D=0$) displayed in (E). Their amplitude is set at $10\%$ of the stationary pulse amplitude. (B) $I_0=3$ and $\tau_D=1.5$: asymmetric slosher: snapshots of the time evolution and periodic evolution of the amplitude at zero (B'). (C) $I_0=3$ and $\tau_D=0.5$, yields a stationary pulse. (D) corresponds to $I_0=0.5$ and $\tau_D=1$ and yields a symmetric breather. Movies of the evolution are provided in supplementary material. Colormap was not thresholded for legibility. }
	\label{fig:Sloshers2d}
\end{figure}

Equating modulus and argument in \eqref{eq:HigherOrder} provide the explicit formulae for the destabilization of the pulse along an asymmetric mode with $n$ oscillations along the unit circle.

The shape of the Hopf bifurcation curves are much less easy to grasp than in the one dimensional case, and we resort to numerical analysis in order to characterize those bifurcations. In Figure~\ref{fig:Bump2d}, we plotted the radially symmetric Hopf bifurcation curve (red line). This curve still connects to the line $I_0=\theta$, but in that case, the delay necessary to induce a destabilization is no more zero, as was the case for the one-dimensional pulse instability. This is visible in our formulae: indeed, while the vanishing of the differential of the pulse at $I_0=\theta$ implied an infinite value for $\omega$, this is no more the case because of the factor $a^{d-1}$ that also vanishes at $I_0=\theta$ and that may compensate the cancellation of $ \U'(a)$. This is indeed what happens in dimension 2. For the parameters chosen, no instability along asymmetric modes were found, and for asymmetric initial conditions and parameters in the instability region, the solutions converge towards the symmetric breathing pulse. 

Of course, other choices of parameters lead to destabilizing the higher-order spatial modes (i.e. destabilizations of the mode $n$ given by~\eqref{eq:HigherOrder}). An example in two dimensions is provided in Fig.~\ref{fig:Sloshers2d}. For the parameters chosen, we find a curve of Hopf bifurcations along the symmetric mode (blue line) together with destabilization of higher order modes. These bifurcation curve decompose the parameter space into 5 regions depending on the stability of the symmetric and asymmetric modes, and asymmetric periodic patterns arise. 

\medskip

We have therefore shown that delays shape the spatio-temporal evolution of the neural fields in the presence of localized input. We now show that this is also the case of fronts. 

%% file: Front.tex
\begin{figure}
	\centering
		\subfigure[Input $\I(x)$.]{\includegraphics[width=.4\textwidth]{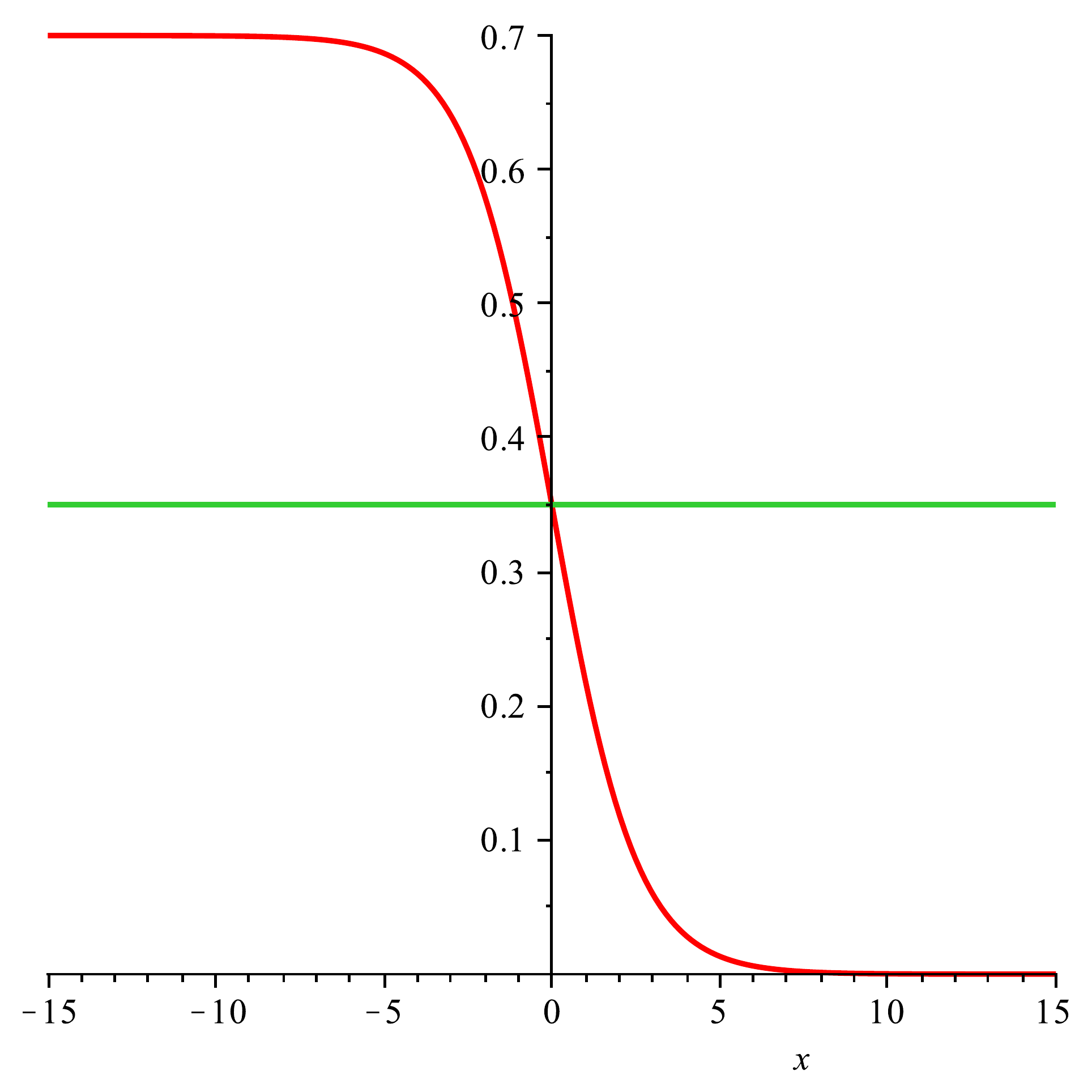}}\qquad\quad\quad
		\subfigure[Inverse Mexican Hat $\J_1(|x|)$.]{\includegraphics[width=.4\textwidth]{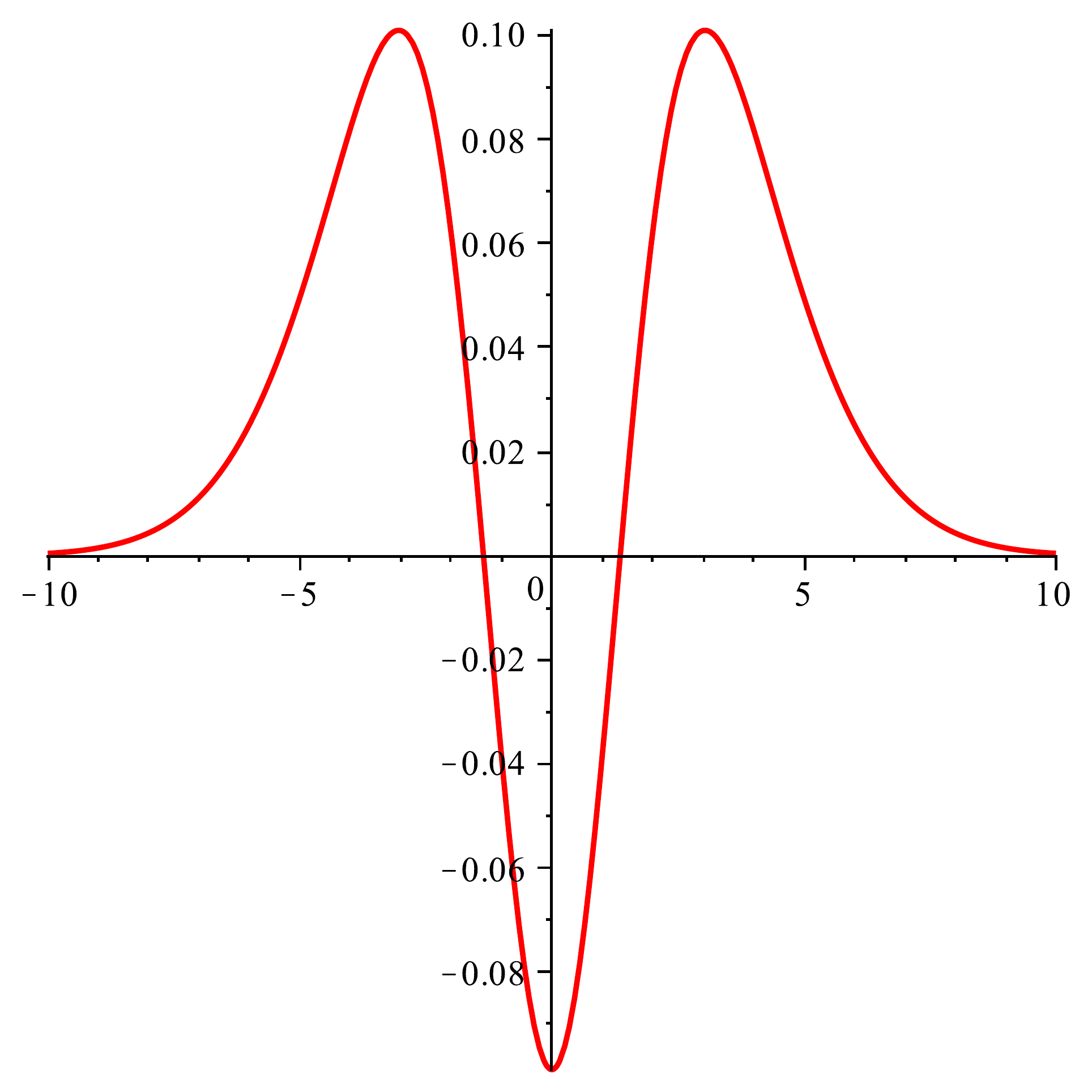}}
	\caption{Form of the input $\I(x)$ (a) and the connectivity function $\J_1(|x|)$ for typical values of the parameters: $I_0=0.7$, $s=0.79$, $w_e=1.7$, $w_i=1.2$, $\sigma_e=4$ and $\sigma_i=2$.}
	\label{fig:InputConnectivityFront1D}
\end{figure}

Let us now consider equation~\eqref{eq:NFE} in dimension $d\geq 1$ with an input of the form
\bqs
I(\mx)=\I(\me\cdot\mx)=I_0\left(1-\frac{1}{1+e^{-s \me \cdot \mx}}\right),
\eqs
for a given direction $\mathbf{e}\in \mathbb{S}^{d-1}$. Here, $I_0$ represents the strength of the input and $s$ the stiffness of the slope of the input at $x=0$ (see Fig.~\ref{fig:InputConnectivityFront1D}(a)). A stationary planar front $\mathcal{V}(x)$ is a time-independent profile of the variable $x=\mathbf{e} \cdot \mathbf{x} $ solution of the equation~\eqref{eq:NFE} satisfying the limits $\underset{x\rightarrow - \infty}{\lim}\mathcal{V}(x)=w_e-w_i+I_0$ and $\underset{x\rightarrow + \infty}{\lim}\mathcal{V}(x)=0$. Throughout this section we suppose that $w_e-w_i+I_0>\theta$. Without loss of generality, we suppose that $\me=(1,0,\dots,0)$ so that $\mathcal{V}(x)$ is solution of
\begin{equation}\label{eq:Front}
\V(x) =  \int_{\R} \J_1(|x-y|)f(\V(y)) dy+\I(x), \quad x\in\R.
\end{equation}
Here, the one-dimensional kernel $\J_1$ appears thanks to our specific scaling of the Gaussian $\G_d$ through which the connectivity $\J_d$ is defined. In particular, this kernel $\J_1$ is defined through the same parameters $(w_e,\sigma_e,w_i,\sigma_i)$ as was $\J_d$. 

\begin{figure}
	\centering
		\subfigure[Profile of $\V(x)$.]{\includegraphics[width=.4\textwidth]{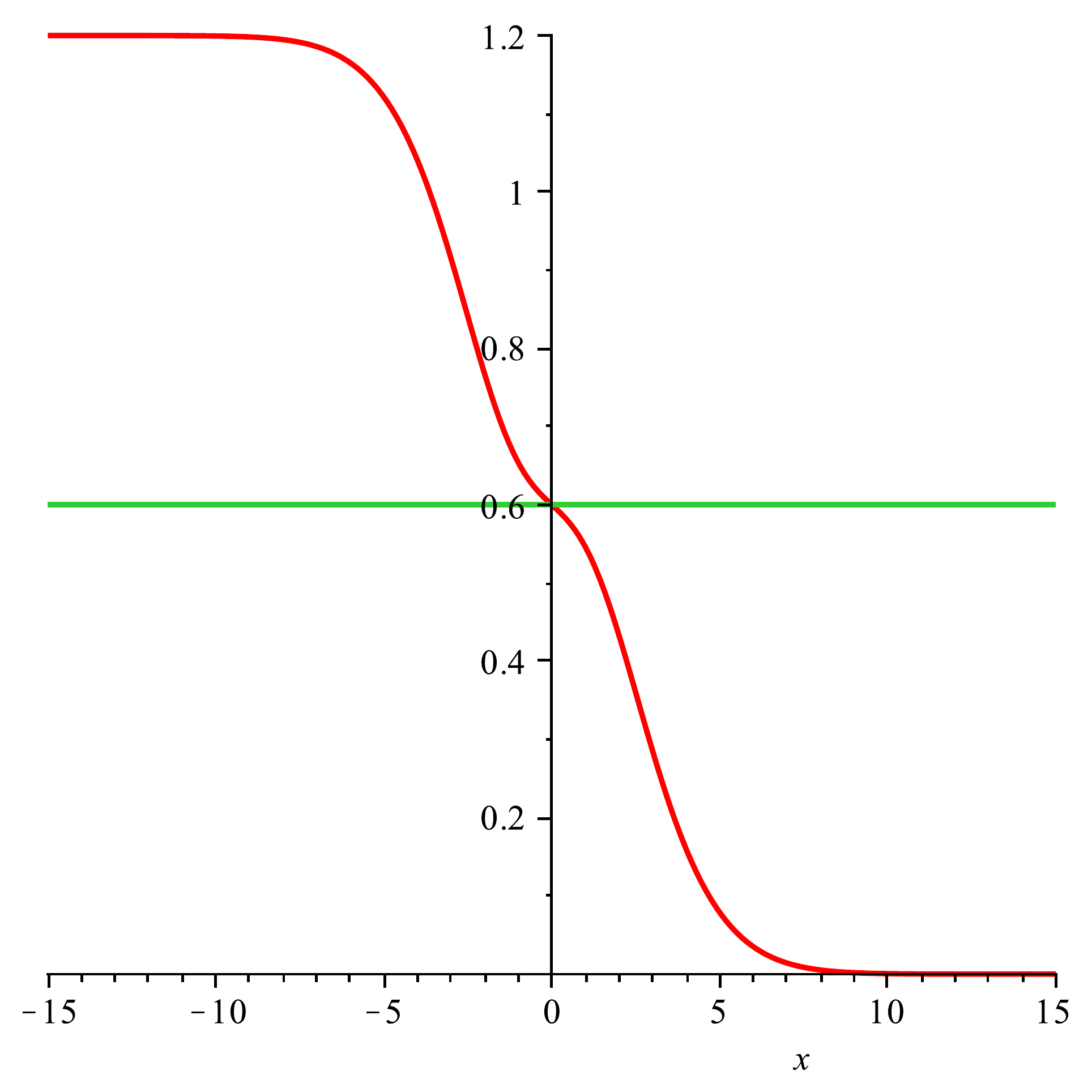}}\qquad\quad\quad
		\subfigure[Zoom around $x=0$.]{\includegraphics[width=.4\textwidth]{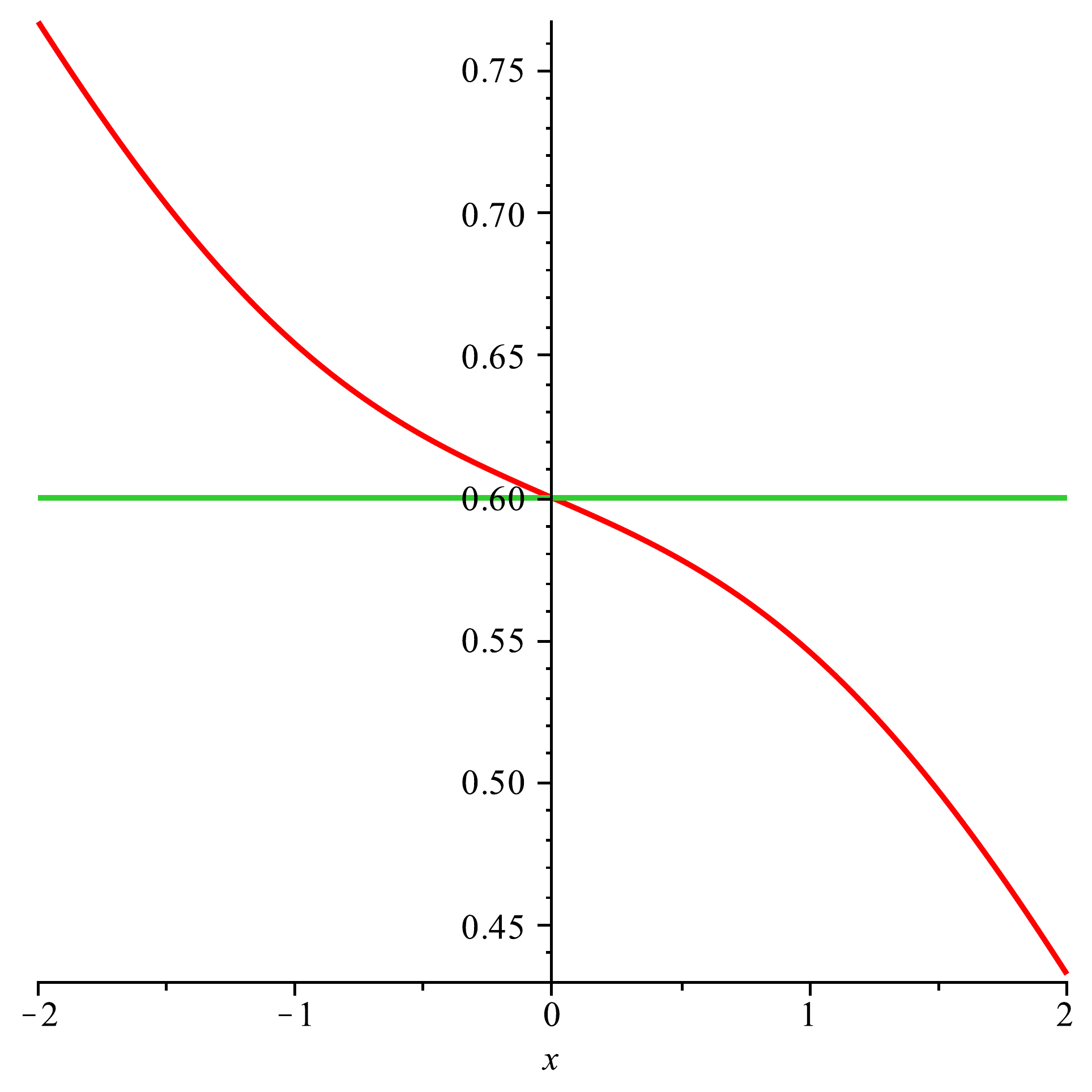}}
	\caption{Profile of the monotonous front for typical values of the parameters: $I_0=0.7$, $s=0.79$, $w_e=1.7$, $w_i=1.2$, $\sigma_e=4$ and $\sigma_i=2$.}
	\label{fig:ShapeFronts1DM}
\end{figure}

To simplify the presentation, we will suppose that planar fronts always cross the threshold $\theta$ at $x=0$, even in the case $I_0=0$ (no input) for which the equation \eqref{eq:Front} is translation invariant\footnote{In the absence of input, there exists a continuum of solutions if one does not specify the location of the crossing value}. Depending on the parameters of the connectivity function $\mathcal{J}_1$, the planar front $\mathcal{V} (x)$ may have different forms. First, it can be monotonous (see Fig.~\ref{fig:ShapeFronts1DM}), in which case we simply have 
\bqq
\label{eq:FM}
\V(x)=\int_x^{+\infty} \J_1(|y|)dy+\I(x),
\eqq
and it crosses the threshold  $\theta$ exactly once. This gives us a relationship between coefficients that needs to be satisfied in order to ensure the existence of such fronts. This relation is given by
\bqq
\label{eq:RTM}
\theta = \frac{w_e-w_i}{2}+\frac{I_0}{2}.
\eqq
Monotonous front occurs only when the Mexican hat function $\J_1$ is greater than the derivative of the input, as $\V'(x)=-\J_1(|x|)+\I'(x)\leq0$ for all $x\in\R$ (see Fig.~\ref{fig:InputConnectivityFront1D}(b)).

\begin{figure}
	\centering
				\subfigure[Profile of $\V(x)$.]{\includegraphics[width=.4\textwidth]{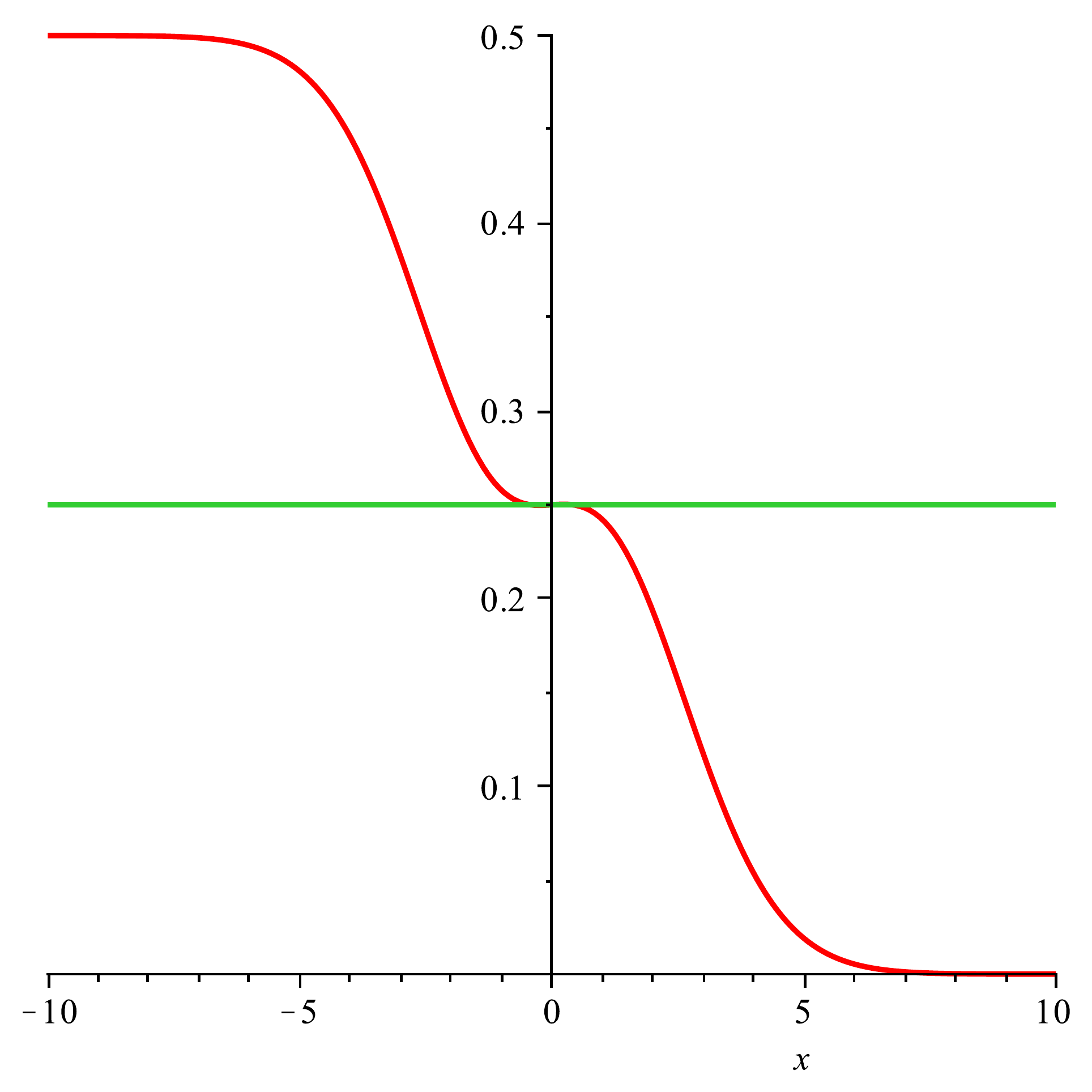}}\qquad\quad\quad
		\subfigure[Zoom around $x=0$.]{\includegraphics[width=.4\textwidth]{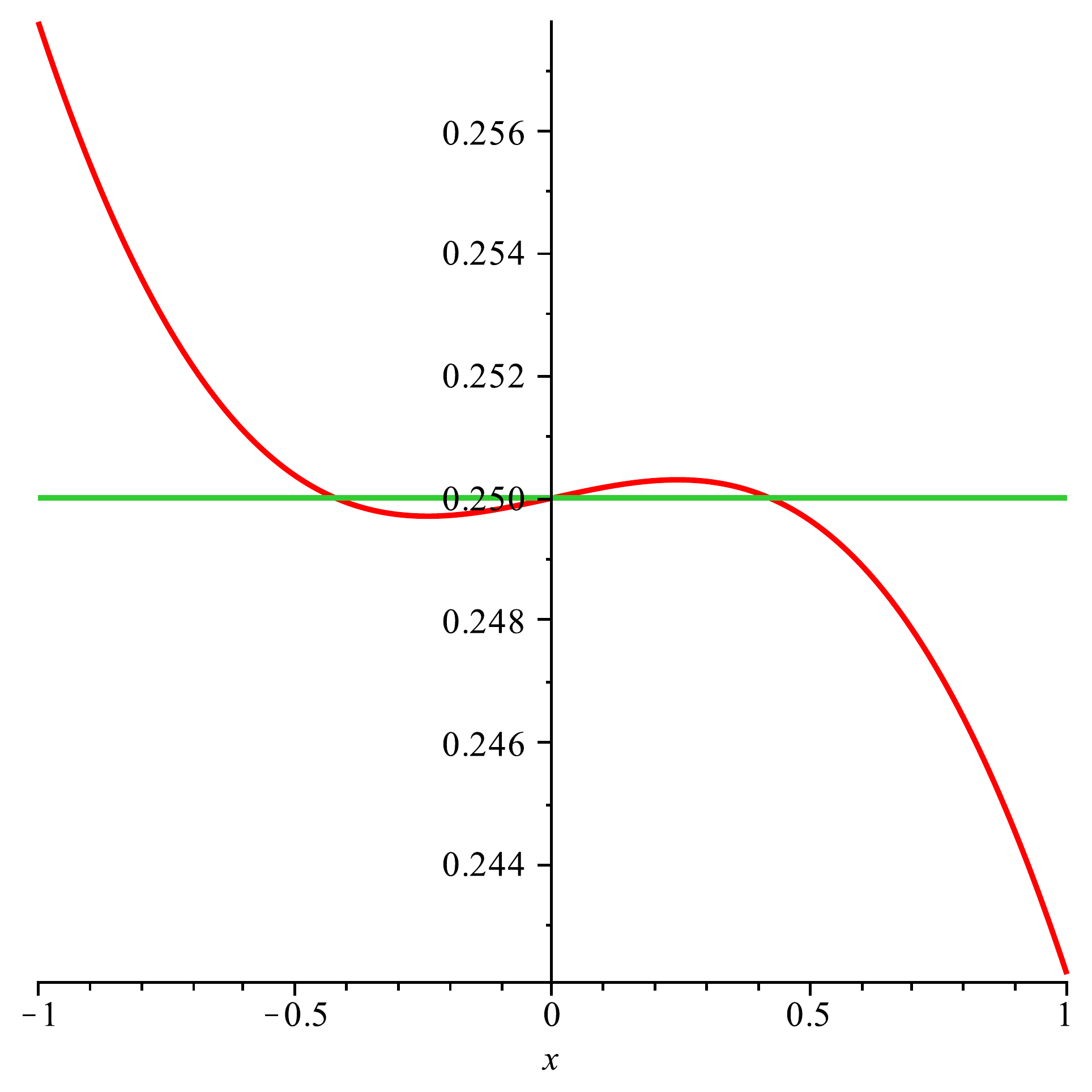}}
	\caption{Profile of the non-monotonous front for typical values of the parameters:  $w_e=1.7$, $w_i=1.2$, $\sigma_e=3$ and $\sigma_i=2$. Here, the value of the crossing is $a\approx 0.422$.}
	\label{fig:ShapeFronts1DNM}
\end{figure}

The second and more difficult case is when the front is non-monotonous, in which case it can cross the threshold $\theta$ an arbitrary number of times. In this paper, we study the specific case when the profile of the front crosses the threshold  exactly three times, at $x=- a$, $x=0$ and $x=a$ for $a>0$ (see Fig.~\ref{fig:ShapeFronts1DNM}). The analysis in the case where the front crosses the threshold more than three times can easily been extended from our analysis. The profile $\V(x)$ of the front now satisfies
\bqs
 \V(x)=\int_{- \infty}^{-a}\J_1(|x-y|) dy + \int_{0}^{a}\J_1(|x-y|) dy +\I(x),
\eqs
or equivalently
\bqq
\label{eq:FNM}
 \V(x)=\int_{x-a}^{x}\J_1(|y|) dy + \int_{x+a}^{\infty}\J_1(|y|) dy +\I(x).
\eqq
Setting $x=\pm a$ and $x=0$ in equation \eqref{eq:FNM}, we respectively find an implicit equation for $a$:
\bqq
\label{eq:RAM}
0=\int_a^{2a} \J_1(y)dy+\frac{\I(-a)-\I(a)}{2},
\eqq
and the same compatibility condition \eqref{eq:RTM} as the one found for monotonous fronts.

\subsection{Linear stability analysis for monotonous fronts}

In order to characterize the linear stability of monotonous front $\V(x)$ with only one crossing at $x=0$, we analyze the spectrum of the linearized operator around the front. Following the analysis of the previous section, it is relatively easy to show that formal differentiation of the Heaviside function $f$ yields to the linearized operator $\cl_d$ around the front solution:
\bqq
\label{eq:OpD}
\cl_d \vp(\mx,t)=-\vp(\mx,t)+\frac{1}{|\V'(0)|}\int_{\R^{d-1}} \J_d(\|\mx-(0,\tilde\my)\|) \varphi((0,\tilde\my),t-\tau( \mx -(0,\tilde\my)))d\tilde\my.
\eqq

In the following sections, we study the spectrum of this operator gaining first insight with the one dimensional case and then turning to more difficult problem of general dimensions.

\subsubsection{One-dimensional neural fields}\label{sec:Mon1D}

In one dimension, the linear operator can be reduced further to
\bqq
\label{eq:OpD1}
\cl_1 \vp(x,t)=-\vp(x,t)+\frac{1}{|\V'(0)|} \J_1(|x|) \varphi(0,t-\tau( x)).
\eqq
Similarly to the case of pulses, we consider complex perturbations in the class of functions $\vp(x,t)=e^{\lambda t}p(x)$, for $\lambda\in\C$ and $p(x)$ a bounded continuous function on $\R$ decaying to zero at infinity, and obtain the dispersion relationship:
\bqq
\label{eq:disp1}
e^{\lambda\tau_D}(\lambda+1)p(x)=\frac{\J_1(|x|)}{|\V'(0)|}p(0)e^{-\lambda\frac{|x|}{c}}.
\eqq
Eigenvalues of the linearized operator are complex solutions $\lambda$ of this transcendental equation.

First, we notice that, once again, the point spectrum is given by the solutions to the equation 
\bqs
\lambda+1=0,
\eqs
and does not contribute to any instability. To find the point spectrum, we estimate the dispersion relation at $x=0$ and obtain the equation
\bqq
\label{eq:disp2}
e^{\lambda\tau_D}(\lambda+1)=\frac{\J_1(0)}{|\V'(0)|}=\frac{\J_1(0)}{\J_1(0)-\I'(0)}.
\eqq
Here, we have used the fact that at $x=0$, $\V'(0)=-\J_1(0)+\I'(0)<0$. A that point, we can already note that relationship \eqref{eq:disp2} is independent of the propagation delay $c$.

\paragraph{\textbf{Case $I_0=0$}} For neural fields equation with no input, the dispersion relation \eqref{eq:disp2} further simplifies to
\bqs
e^{\lambda\tau_D}(\lambda+1)=1.
\eqs
Writing $\lambda=\alpha+i\omega$ and taking the modulus of both sides of this equation, we obtain the relation
\bqs
e^{\alpha \tau_D} \sqrt{(1+\alpha)^2+\omega^2}=1,
\eqs
which stands only if $\alpha \leq 0$ since $\tau_D\geq 0$. Checking that the case $\alpha=0$ implies $\omega=0$ (reflecting the translational invariance of the system), this shows that monotonous stationary fronts are always spectrally stable in the absence of stimulus and neither constant nor propagation delays can destabilize the system.

\begin{figure}
	\centering
		\includegraphics[width=.8\textwidth]{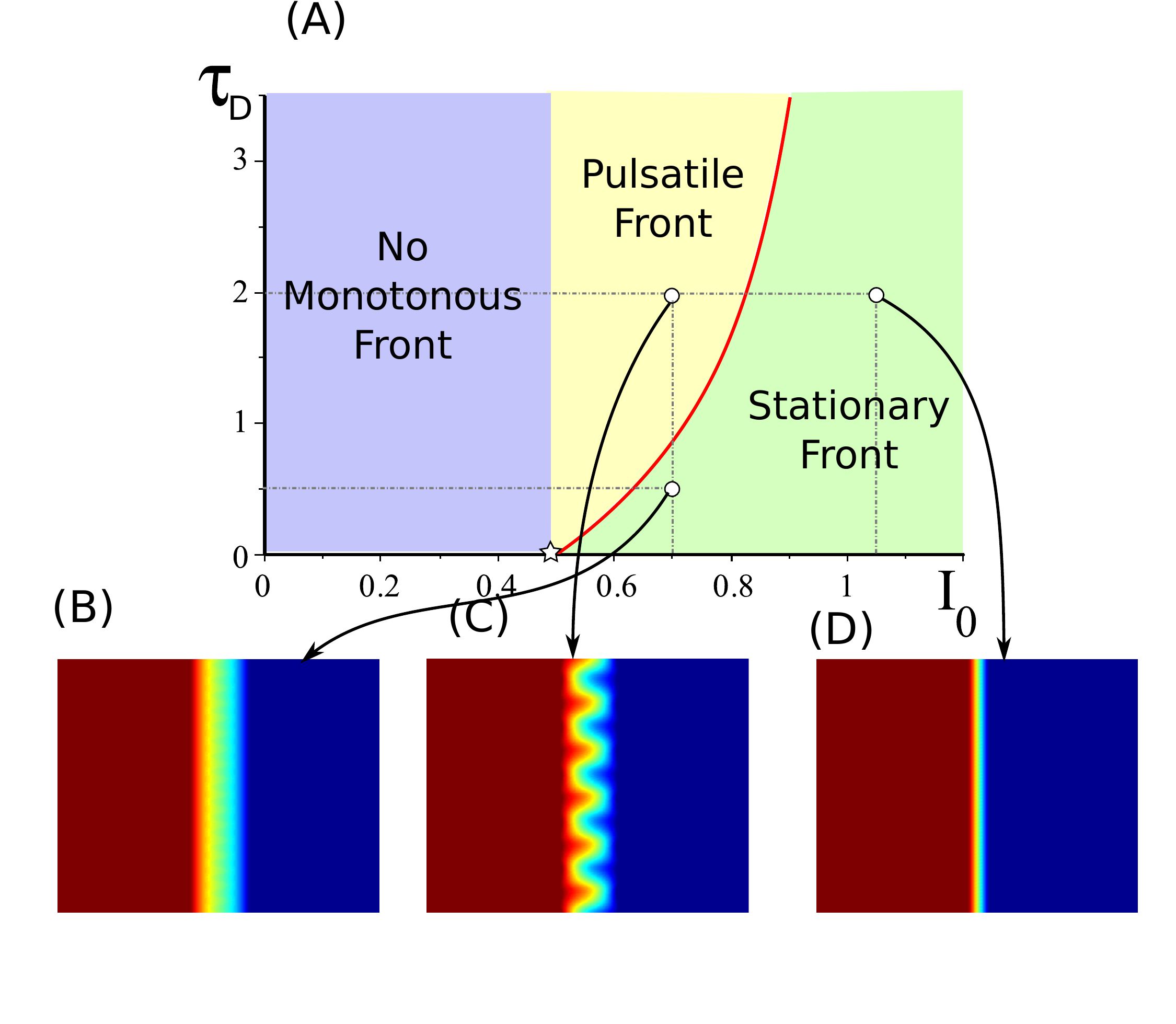}
	\caption{Instabilities of the fronts in the presence of input. (A): bifurcation diagram as a function of the input amplitude $I_0$ and delay $\tau_D$ for an inverse Mexican hat $w_e=1.7$, $w_i=1.2$, $\sigma_e=4$, $\sigma_i=2$. The value of $s$ is taken to be $-8\J_1(0)$. The red line represents the Hopf bifurcation curve, (B) $I_0=0.7$, $\tau_D=0.5$: stationary pulse, (C): $I_0=0.7$, $\tau_D=2$: sloshing instability, (D): $I_0=1.2$, $\tau_D=2$: stationary pulse. Color bar within $\theta \pm 10\%$.  }
	\label{fig:Fronts1dMonotonous}
\end{figure}

\paragraph{\textbf{Case $I_0>0$}} When the input $\I(x)$ is turned on, it is possible to destabilize stationary front solutions through Hopf bifurcations by increasing the constant delay $\tau_D$. From our dispersion relationship \eqref{eq:disp2}, we first see that this can only occur when the connectivity function satisfies $\J_1(0)<0$. In our numerical example, we will chose an inverse Mexican hat which naturally satisfies this condition. In that case, Hopf bifurcations arise if:
\bqq
\label{eq:CondHopf1}
\frac{\J_1(0)}{\J_1(0)-\I'(0)}<-1,
\eqq
and for parameters satisfying the relationship
\begin{equation}\label{eq:DelayFront1}
\left\{ 
\begin{array}{l}
\omega = \sqrt{ \left(\frac{\J_1(0)}{\J_1(0)-\I'(0)}\right)^2  -1}.  \\ 
{\tau_{D}} = \frac{1}{\omega} \left ( \pi -\text{arccos}\left(\frac{\J_1(0)-\I'(0)}{|J_1(0)|} \right) \right ).
\end{array} 
\right.
\end{equation}
The condition \eqref{eq:CondHopf1} together with the fact that $\V'(0)<0$ give us a range of allowed values for $I_0$ as a function of the other parameters for which we can find Hopf instabilities. More precisely we have that
\bqs
I_0 \in \left(I_m,I_M \right):=\left(-\frac{4\J_1(0)}{s},-\frac{8\J_1(0)}{s} \right).
\eqs
In particular, this informs us that the destabilization curve given by \eqref{eq:DelayFront1} (see Fig~\ref{fig:Fronts1dMonotonous}(A)) as a function of $I_0$ (while all other parameters are fixed) emerges from $\tau_D=0$ at a value $I_0=I_m$. At the other boundary of this interval, the corresponding value of $\omega$ tends to zero, and therefore the delay $\tau_D$ corresponding to the Hopf bifurcation diverges to infinity. The curve $I_0 \longmapsto \tau_D(I_0)$ given by \eqref{eq:DelayFront1} naturally partitions the parameter space $(I_0,\tau_D)$ into $3$ regions: no monotonous fronts, stationary monotonous fronts and pulsatile fronts. Note that the mode associated to the Hopf instability is always symmetric (see \eqref{eq:disp1}) and that the notion of \textit{breathers} or \textit{sloshers} is not appropriate for monotonous fronts. This is why we use the terminology \textit{pulsatile} fronts and panel $(C)$ of Fig~\ref{fig:Fronts1dMonotonous} shows such a periodic front.

\subsubsection{General neural fields}

We can now extend the previous analysis in the case of general neural fields as follows in the absence of external input. First, in order to obtain closed-form formula, we only consider the case of neural fields with constant delays $\tau_D$. In that case, the neural field equation is simply given by
\bqq
\label{eq:nfe2d}
\frac{\partial u}{\partial t}(\mx,t) = -u(\mx,t) + \int_{\R^d} \J_d(\|\mx-\mx'\|)f(u(\mx',t-\tau_D)) d\mx',
\eqq
for all $\mx=(x,\tilde{\mx})\in\R\times\R^{d-1}$. As seen in the introduction of the section, we suppose that there exists a stationary planar solution $u(\mx):=\V(x)$ of \eqref{eq:nfe2d}, with profile $\V$ solution of 
\bqs
\V(x) = \int_\R \J_1(|x-x'|)f(\V(x'))dx',
\eqs
where we again used the fact that $\J_1(x)=\int_{\R^{d-1}} \J_{d}(\| (x,\tilde{\mx}) \|)d\tilde{\mx}$. In order to ensure the existence of a monotonous profile $\V(x)$, we suppose that $\J_1(|x|)\geq0$ for all $x\in\R$ while $\theta$ is related to other parameters via \eqref{eq:RTM} with $I_0=0$.

\begin{figure}
	\centering
\subfigure[Profile of $\J_2(\|\mx\|)$.]{\includegraphics[width=.5\textwidth]{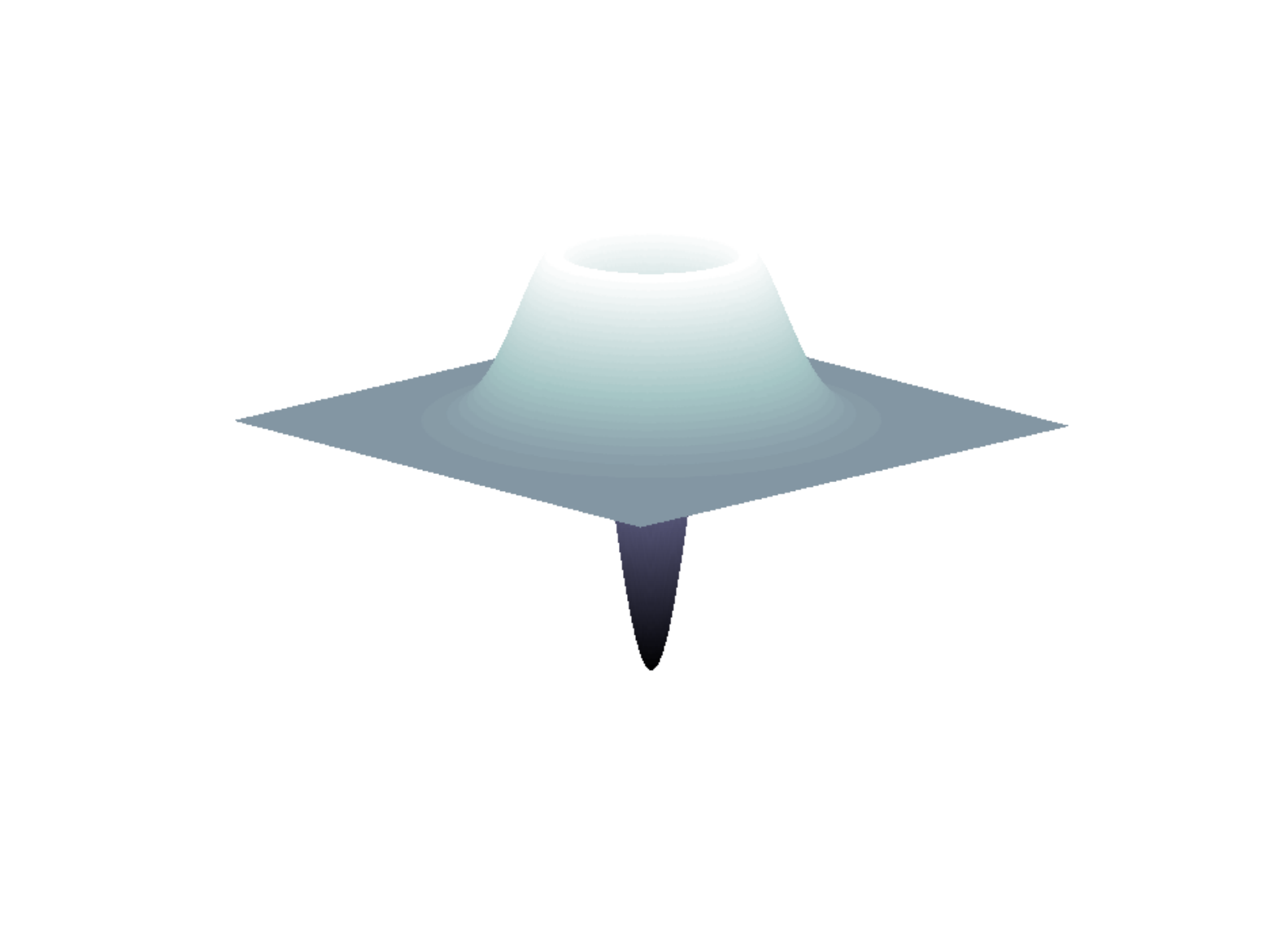}}\qquad\quad\quad
		\subfigure[Profile of $\J_1(|x|)$.]{\includegraphics[width=.3\textwidth]{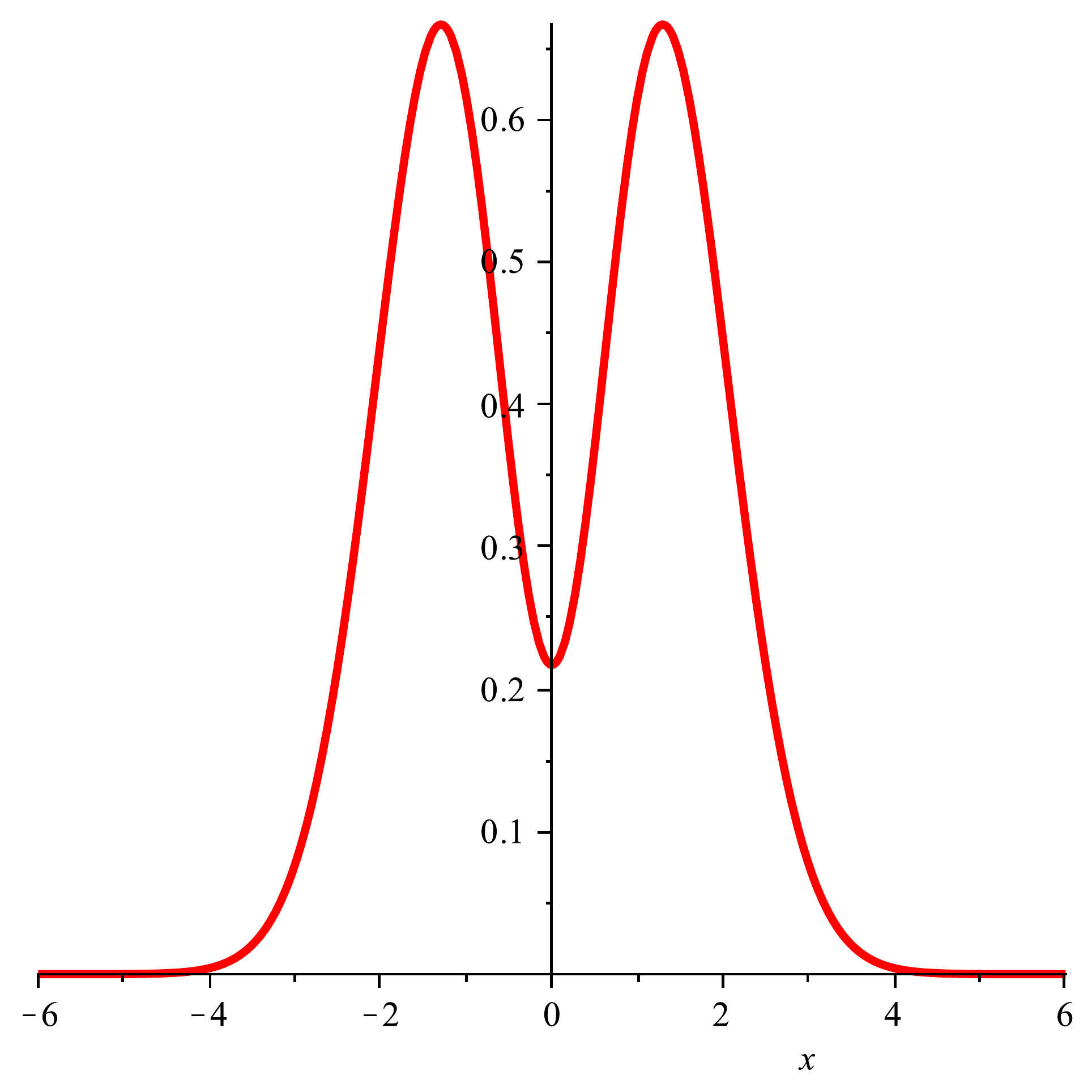}}
	\caption{Profiles of $\J_2$ and $\J_1$ for values of the parameters given by  $w_e=15$, $w_i=12.5$, $\sigma_e=1.5$ and $\sigma_i=1.3$.}
	\label{fig:ShapeConnectivities}
\end{figure}

The linear operator given in \eqref{eq:OpD} simplifies to
\bqs
\cl_d \vp(\mx,t)=-\vp(\mx,t)+\frac{1}{|\V'(0)|}\int_{\R^{d-1}} \J_d(\|\mx-(0,\tilde\my)\|) \varphi((0,\tilde\my),t-\tau_D))d\tilde\my.
\eqs
We now consider complex perturbations in the class of functions $\vp(\mx,t)=e^{\lambda t+i \ell\cdot \tilde\mx}p(x)$, where $\mx=(x,\tilde\mx)$ and $\ell\in\R^{d-1}$, for $\lambda\in\C$ and $p(x)$ a bounded continuous function on $\R$ decaying to zero at infinity. On this class of functions, we obtain the dispersion relationship:
\bqq
\label{eq:disp3}
e^{\lambda\tau_D}(\lambda+1)p(x)=\frac{1}{|\V'(0)|}p(0)\left[w_e\G_1(x,\sigma_e)\widehat{\G}_{d-1}(\ell,\sigma_e)-w_i\G_1(x,\sigma_i)\widehat{\G}_{d-1}(\ell,\sigma_i) \right], 
\eqq
for $\ell\in\R^{d-1}$, where $\widehat\G_{d-1}$ is defined as 
\bqs
\widehat\G_{d-1}(\ell,\sigma)= \exp\left(-\frac{\sigma^2\|\ell\|^2}{4}\right), \quad \forall \ell\in\R^d.
\eqs
For non zero perturbations $p$ that satisfy $p(0)=0$, we have that $\lambda=-1$, which gives the essential spectrum. For perturbations with $p(0)\neq0$, we are left to solve the equation
\bqq
\label{eq:disp4}
e^{\lambda\tau_D}(\lambda+1)= \frac{1}{|\V'(0)|}\left[w_e\G_1(0,\sigma_e)\widehat{\G}_{d-1}(\ell,\sigma_e)-w_i\G_1(0,\sigma_i)\widehat{\G}_{d-1}(\ell,\sigma_i) \right]:=\Psi(\|\ell\|).
\eqq
We observe that $\Psi(0)=1$ and $\Psi(x)\rightarrow 0$ as $x \rightarrow +\infty$. One can also easily check by direct computations that $\J_1(0)>0$ implies $\Psi(\|\ell\|)\leq1$ for all $\ell \in \R^{d-1}$ and that $0$ is the only maximum of $\Psi$. 

We are interested in changes in the stability of the planar front solutions. Fold transitions arise under the condition that
\bqs
\Psi(\|\ell\|)=1.
\eqs
This condition is satisfied only at $\ell=0_{\R^{d-1}}$ and we recover the translation invariance in the direction of the front. We know look for destabilization of the planar front through Hopf bifurcations. A necessary and sufficient condition for the rightmost characteristic eigenvalue of \eqref{eq:disp4} to be purely imaginary and nonzero is the existence of $l_0>0$ such that
\bqs
\Psi(l_0)<-1.
\eqs
In this case the corresponding frequency $\omega$ and associated critical constant delay $\tau_D$ satisfy the relations
\begin{equation}\label{eq:DF}
\left\{ 
\begin{array}{l}
\omega = \sqrt{  \Psi(l_0)^2  -1}.  \\ 
{\tau_{D}} = \frac{1}{\omega} \left ( \pi -\text{arccos}\left(\frac{1}{|\Psi(l_0)|} \right) \right ).
\end{array} 
\right.
\end{equation}
Interestingly, this loss of stability due to delays arises only in higher dimensions and we obtain a full sphere of critical points $l_0=\|\ell\|$, $\ell\in\R^{d-1}$, due to the rotational symmetry in the transverse directions. 

\begin{figure}[htbp]
	\centering
		\includegraphics[width=1\textwidth]{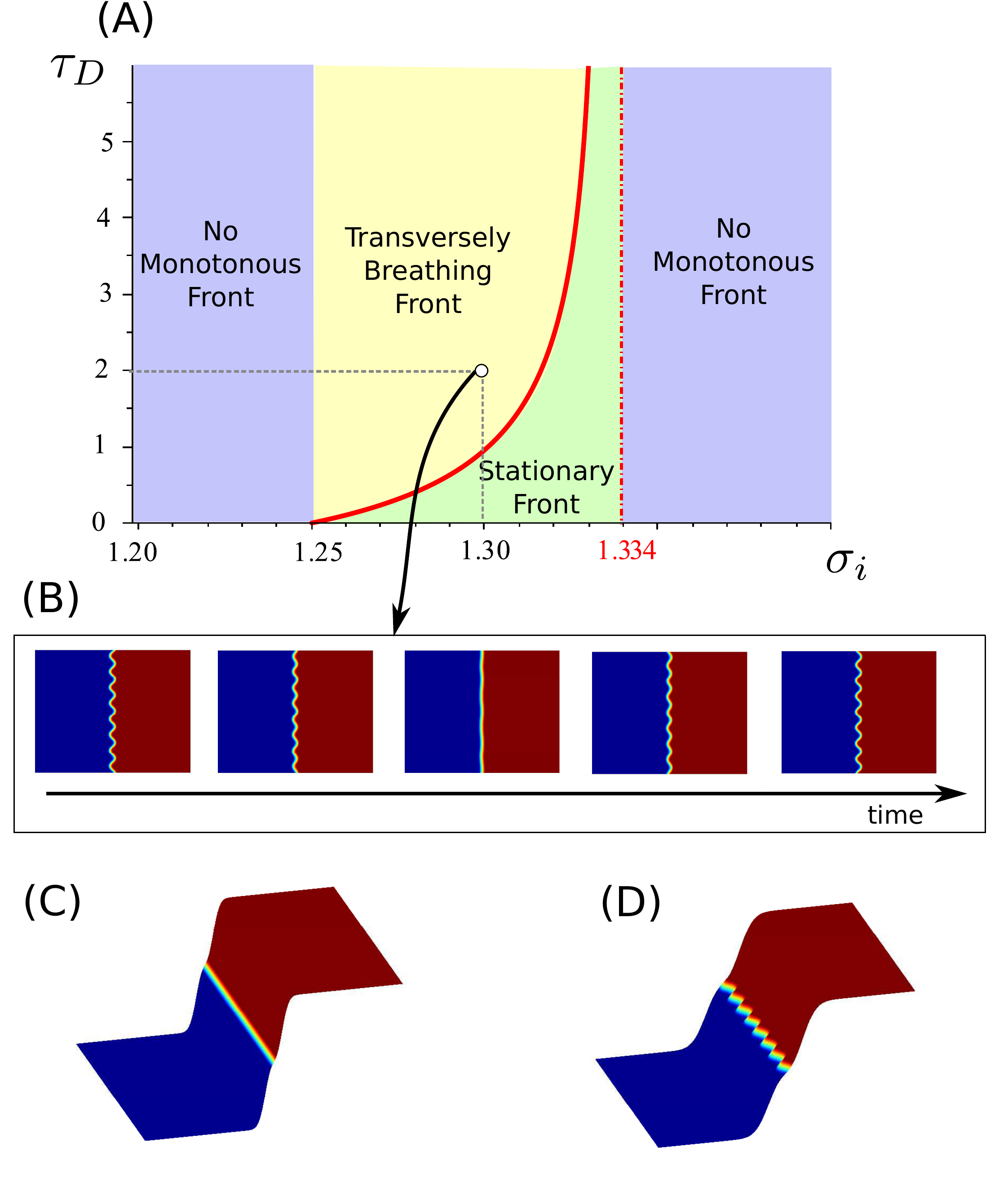}
	\caption{Stability of planar stationary fronts for an inverse Mexican hat connectivity: $w_e=15$, $w_i=12.5$, $\sigma_e=1.5$ as constant delay $\tau_D$ and $\sigma_i$ are varied. (A) Bifurcation diagram in the $(\sigma_i,\tau_D)$ plane with the first Hopf bifurcation curves given by equations \eqref{eq:DF} (red curve). This curve emerges at $\sigma_i=1.25$ ($\tau_D=0$) and asymptotes the line $\sigma_i\approx 1.334$ ($\tau_d\rightarrow \infty$). Blue region: non monotonous front, Green: stationary planar front (C) and Yellow: transversely breathing front (D). (B) Five snapshots over half a period of the transversely pulsatile planar front at $\sigma_i=1.3$ and $\tau_D=2$, seen from above. Movies of the evolution are provided in supplementary material. }
	\label{fig:Transverse}
\end{figure}

We illustrate this analysis in the case of $d=2$ with values of the parameters chosen so that $\J_2$ is an inverse Mexican hat connectivity function and $\J_1$ is positive for all values (see Fig.~\ref{fig:ShapeConnectivities}). More precisely, we fix the values of $w_e$, $w_i$, $\sigma_e$ and treat $\sigma_i$ as a bifurcation parameter. We obtain a bifurcation diagram in the $(\sigma_i,\tau_D)$ plane with the first Hopf bifurcation curve given by equations \eqref{eq:DF} as shown in panel (A) of Fig.~\ref{fig:Transverse}. This curve emerges at $(\sigma_i,\tau_D)=(1.25,0)$ and asymptotes the vertical line $\sigma_i\approx 1.334$ where $\tau_D\rightarrow +\infty$. There exists a region in parameter space where stationary monotonous planar fronts exist and can be destabilized by increasing the value of the constant delay $\tau_D$. In that case, transverse modes lead to the destabilization of the planar front through a Hopf bifurcation and a sequence of snapshots (see Fig.~\ref{fig:Transverse}(B,D)) over half a period illustrates this periodic behavior. This feature was not possible in the one dimensional case, as in the absence of stimulus, stationary fronts are always spectrally stable (see Sec.\ref{sec:Mon1D}). This clearly indicates that dimension of the neural fields and delays play a crucial role in the emergence of structured spatio-temporal cortical activities. This also suggests that the information obtained in one dimensional models, seen as a simplification of higher dimensional systems, can sometimes be misleading.

\subsection{Linear stability analysis for non-monotonous fronts}

In this section, we study the linear stability of non-monotonous fronts $\V(x)$ crossing the threshold $\theta$ exactly three times at $x=\pm a $ and $x=0$ in the absence of stimulus. In that case, $a$ is simply defined through equation 
\bqs
0=\int_a^{2a}\J_1(x)dx.
\eqs
We denote $\E=\{\pm a,0\}$. For non-monotonous fronts, the linear operator $\cl_d$ is modified according to
\bqq
\label{eq:OpDNM}
\cl_d \vp(\mx,t)=-\vp(\mx,t)+\sum_{\xi \in \E }\frac{1}{|\V'(\xi)|}\int_{\R^{d-1}} \J_d(\|\mx-(\xi,\tilde\my)\|) \varphi((\xi,\tilde\my),t-\tau(\mx -(\xi,\tilde\my) ))d\tilde\my.
\eqq

In order to simplify our presentation, we only consider the case of one-dimensional non-monotonous fronts. As a consequence, in one dimension, $\cl_1$ further reduces to
\bqq
\label{eq:OpDNM1}
\cl_1 \vp(x,t)=-\vp(x,t)+\sum_{\xi \in \E}\frac{1}{|\V'(\xi)|} \J_1(|x-\xi|) \varphi(\xi,t-\tau( x-\xi)).
\eqq
As usual, we are interested in studying perturbations of the form $\vp(x,t)=e^{\lambda t}p(x)$, for $\lambda\in\C$ and $p(x)$ a bounded continuous function on $\R$ decaying to zero at infinity. Replacing this Ansatz into $\partial_t \vp =\cl_1 \vp$ we obtain
\bqq
\label{eq:disp5}
e^{\lambda \tau_D}(\lambda+1)p(x)=\frac{\J_1(x)}{|\V'(0)|}p(0)e^{-\lambda\frac{|x|}{c}}+\frac{\J_1(|x-a|)}{|\V'(a)|}p(a)e^{-\lambda\frac{|x-a|}{c}}+\frac{\J_1(|x+a|)}{|\V'(-a)|}p(-a)e^{-\lambda\frac{|x+a|}{c}},
\eqq
Estimating equation \eqref{eq:disp5} at $x=\pm a$ and $x=0$ yields the system
\bqq
\label{eq:Syst1}
\left( 
\begin{matrix}
\mG(\lambda)-\gamma_1\J_1(0) & -\gamma_2 \J_1(a)e^{-\frac{\lambda a}{c}} & -\gamma_1\J_1(2a)e^{-\frac{2 \lambda a}{c}}  \\ -\gamma_1 \J_1(a) e^{-\frac{\lambda a}{c}} & \mG(\lambda)-\gamma_2\J_1(0) & -\gamma_1\J_1(a) e^{-\frac{\lambda a}{c}}
 \\ -\gamma_1\J_1(2a)e^{-\frac{2 \lambda a}{c}} & -\gamma_2\J_1(a) e^{-\frac{\lambda a}{c}} & \mG(\lambda)-\gamma_1\J_1(0)
\end{matrix}
\right)
\left( \begin{matrix} p(-a) \\ p(0) \\ p(a)\end{matrix}\right)=0
\eqq
where we defined $\mG(\lambda) := e^{\lambda \tau_D} (\lambda +1)$ and
\begin{align*}
\gamma_1&:=\frac{1}{\J_1(0)+\J_1(2a)-\J_1(a)}>0,\\ 
\gamma_2&:=\frac{1}{\J_1(0)-2\J_1(a)}>0.
\end{align*} 
This matrix system has the form
\bqs
\mathcal{M}=\left(\begin{matrix} m_1 & m_2 & m_3 \\ m_4 & m_5 & m_4 \\ m_3 & m_2 & m_1 \end{matrix} \right),
\eqs
so that its determinant can be factored as $\det(\mathcal{M})=(m_1-m_3)(m_1m_5-2m_2m_4+m_3m_5)$. System \eqref{eq:Syst1} has non trivial solutions if and only if $\det(\mathcal{M})=0$, which leads to two relationships:
\bqq
\label{eq:Relation1}
\mG(\lambda)=\gamma_1\left(\J_1(0)-\J_1(2a)e^{-\frac{2\lambda a}{c}} \right),
\eqq
and
\bqq
\label{eq:Relation2}
\mG(\lambda)^2-\alpha(a,c)\mG(\lambda)+\beta(a,c)=0,
\eqq
where
\begin{align*}
\alpha(a,c)&:=(\gamma_1+\gamma_2)\J_1(0)+\gamma_1\J_1(2a)e^{-\frac{2\lambda a}{c}},\\
\beta(a,c)&:=\gamma_1\gamma_2\left[\J_1(0)^2+\left(\J_1(0)\J_1(2a)-2\J_1(a)^2\right)e^{-\frac{2\lambda a}{c}} \right].
\end{align*}

\begin{figure}
	\centering
		\subfigure[Asymmetric Perturbation]{\includegraphics[width=.4\textwidth]{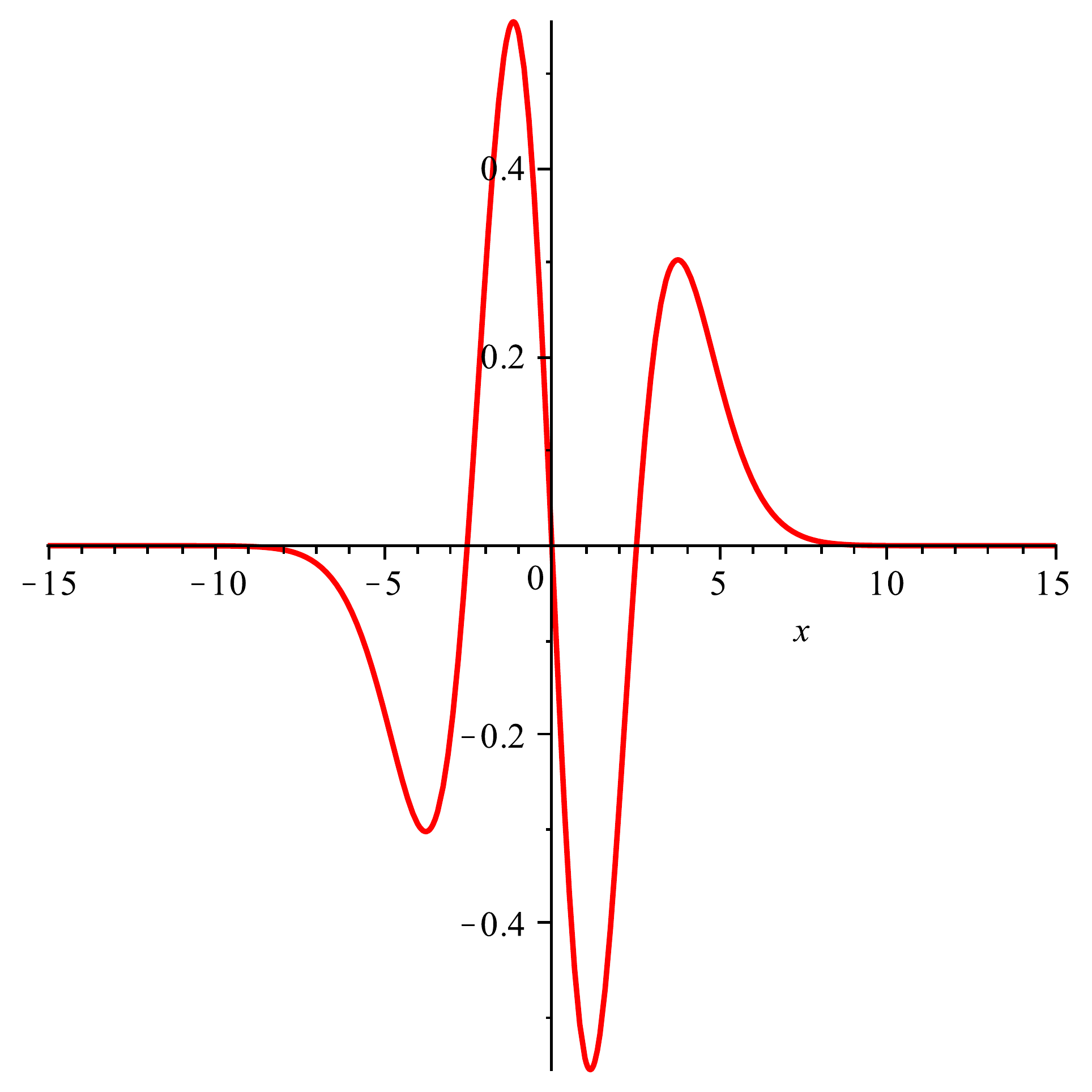}}\qquad\quad\quad
		\subfigure[Symmetric Perturbation]{\includegraphics[width=.4\textwidth]{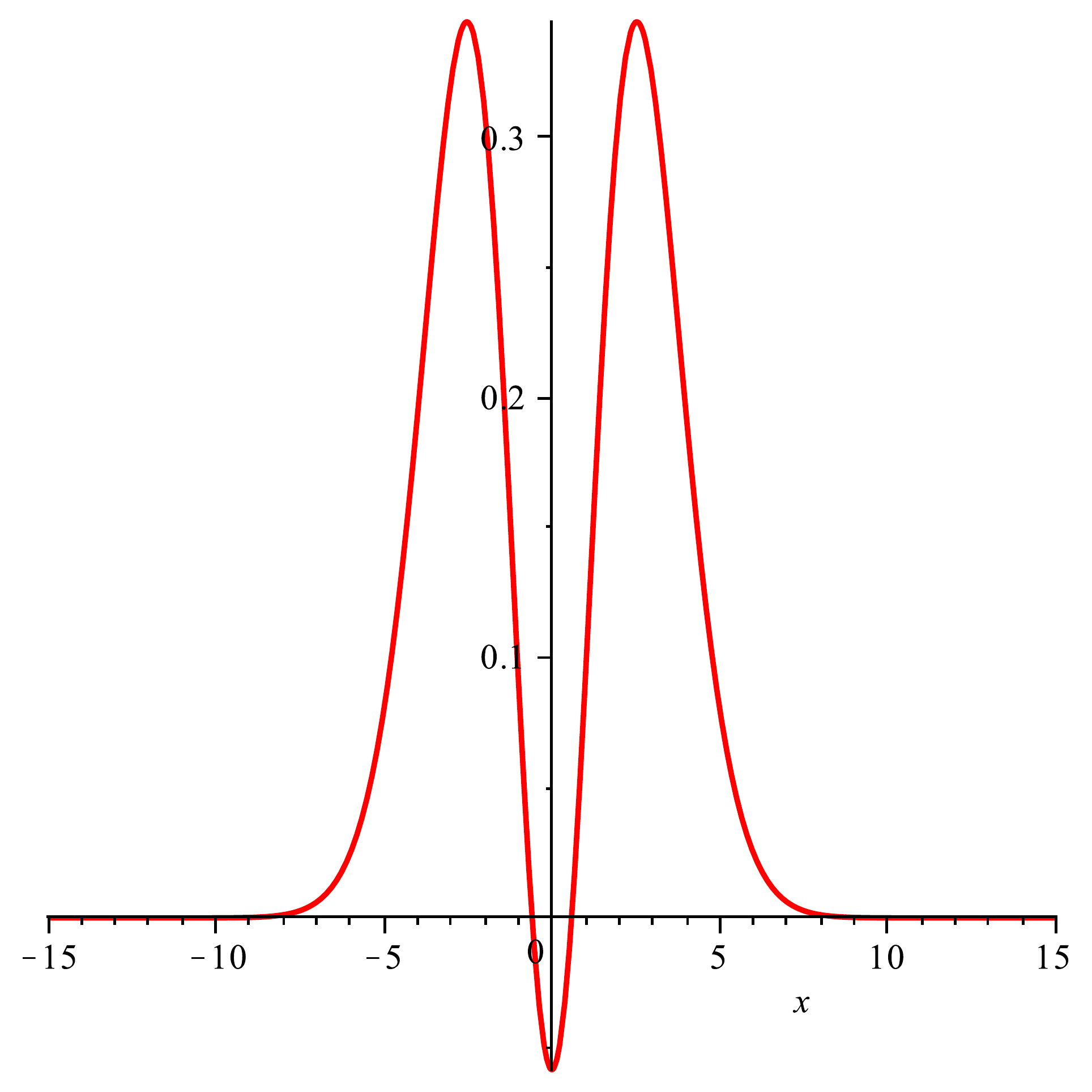}}
	\caption{Form of the eigenfunctions related to the front destabilization for the inverse Mexican hat ($w_e=1.7$, $w_i=1.2$, $\sigma_e=3$ and $\sigma_i=2$): (a) asymmetric perturbation $p_{\text{asym}}$ corresponding to sloshing fronts and (b) symmetric perturbations $p_{\text{sym}}$ corresponding to breathing fronts.}
	\label{fig:PerturbationsFront1D}
\end{figure}

Modes associated to the first relationship \eqref{eq:Relation1} are proportional to
\bqq
\label{eq:Mode1}
p_{\text{asym}}(x)=\J_1(|x-a|)e^{-\lambda\frac{|x-a|}{c}}-\J_1(|x+a|)e^{-\lambda\frac{|x+a|}{c}},
\eqq
and are asymmetric perturbations, while the modes associated to the second relationship \eqref{eq:Relation2} have the form
\bqq
\label{eq:Mode2}
p_{\text{sym}}(x)=\J_1(|x-a|)e^{-\lambda\frac{|x-a|}{c}}+\J_1(|x+a|)e^{-\lambda\frac{|x+a|}{c}}+2\J_1(|x|)e^{-\lambda\frac{|x|}{c}},
\eqq
and are symmetric perturbations. These two eigenfunctions are plotted in Fig.~\ref{fig:PerturbationsFront1D} for a specific choice of parameter. 

First, we study the relationship \eqref{eq:Relation1} and look for possible Hopf bifurcations. In the purely synaptic delay case ($c\rightarrow + \infty$), we obtain the closed relationships:
\begin{equation}
\label{eq:SolRelation1}
\left\{ 
\begin{array}{l}
\omega = \sqrt{\left ( \gamma_1(\J_1(0)-\J_1(2a)) \right ) ^2 -1} \\
\tau_D = \frac{1}{\omega}\left(\pi-\text{arccos}\left(\frac{1}{\gamma_1|\J_1(0)-\J_1(2a)|}\right) \right) , \\
\end{array}
\right.
\end{equation}
provided that
\bqs
\gamma_1(\J_1(0)-\J_1(2a))<-1,
\eqs
or equivalently
\bqs
2\J_1(0)<\J_1(a).
\eqs
In the propagation delay case, we wish again to draw bifurcation curves in the $(\tau_D,c)$ plane. We have: 
\begin{equation}
\label{eq:Relation1bis}
e^{- \frac{2 i a \omega}{c}} = \frac{\J_1(0)-\gamma_1^{-1}e^{i \omega \tau_D} (i \omega +1)}{\J_1(2a)} := \mathcal{Y} (\tau_D,\omega)
\end{equation}
Equating modulus and argument gives :
\begin{equation}
\label{eq:SolRelation1bis}
\left\{ 
\begin{array}{l}
\left | \mathcal{Y} (\tau_D,\omega) \right | = 1 \\
c =  \dfrac{2 a \omega}{2n\pi-\arg(\mathcal{Y} (\tau_D,\omega))}, \quad n\in \N. \\
\end{array}
\right.
\end{equation}

We now return to relationship \eqref{eq:Relation2}. In the purely synaptic case, this relationship is in fact a second order polynomial equation for $\mG(\lambda)$ of the form:
\begin{equation}
\label{eq:Relation2bis}
\mG(\lambda)^2-\alpha_\infty(a) \mG(\lambda)+\beta_\infty(a) = 0
\end{equation}
where 
\bqs
\alpha_\infty(a) = \underset{c\rightarrow \infty}{\lim}\alpha(a,c), \quad \text{ and } \quad \beta_\infty(a) = \underset{c\rightarrow \infty}{\lim}\beta(a,c).
\eqs

Complex numbers $\lambda_\pm$ solutions of \eqref{eq:Relation2bis} thus satisfy: 
\bqs
\mG(\lambda_\pm) = \frac{ \alpha_\infty(a) \pm \sqrt{\alpha_\infty(a)^2-4\beta_\infty(a)}}{2}:=\Theta_\pm(a).
\eqs
There again, we write $\lambda_\pm = i \omega_\pm$ and then equate modulus and argument to study possible Hopf bifurcations. We obtain the two following solutions for $(\omega_\pm,{\tau_{D}}_\pm)$: 

\begin{equation}
\label{eq:SolRelation2bis}
\left\{ 
\begin{array}{l}
\omega_\pm = \sqrt{ \left| \Theta_\pm(a) \right| ^2 -1} \\
{\tau_{D}}_\pm \equiv \frac{1}{\omega_\pm}\left(  \arg\left(\Theta_\pm(a)\right)-\text{arccos}\left( \dfrac{1}{\left| \Theta_\pm(a) \right|}\right)\right) \pmod{\frac{2\pi}{\omega_\pm}}, \\
\end{array}
\right.
\end{equation}
provided that the following condition is satisfied:
\bqs
\left| \Theta_\pm(a) \right| >1.
\eqs
In the propagation delay case, we isolate the term $e^{-2\lambda a / c}$ in equation \eqref{eq:Relation2} to obtain the formula:
\bqs
e^{-2\lambda a / c}= \frac{\mG(\lambda)^2-(\gamma_1+\gamma_2)\J_1(0)\mG(\lambda)+\gamma_1\gamma_2\J_1(0)^2}{\gamma_1\J_1(2a)\mG(\lambda)+\gamma_1\gamma_2 \left(2\J_1(a)^2-\J_1(0)\J_1(2a) \right)}:=\mathcal{Z} (\tau_D,\lambda).
\eqs
Taking the modulus and argument of the above equation for $\lambda=i\omega$ yields the implicit relationships for $c$ and $\tau_D$
\begin{equation}
\label{eq:SolRelation2bisProp}
\left\{ 
\begin{array}{l}
\left | \mathcal{Z} (\tau_D,i\omega) \right | = 1, \\
c =  \dfrac{2 a \omega}{2n\pi- \arg(\mathcal{Z} (\tau_D,i\omega))}.
\end{array}
\right.
\end{equation}
for some $n \in \N$.

\begin{figure}[htbp]
	\centering
		\includegraphics[width=0.8\textwidth]{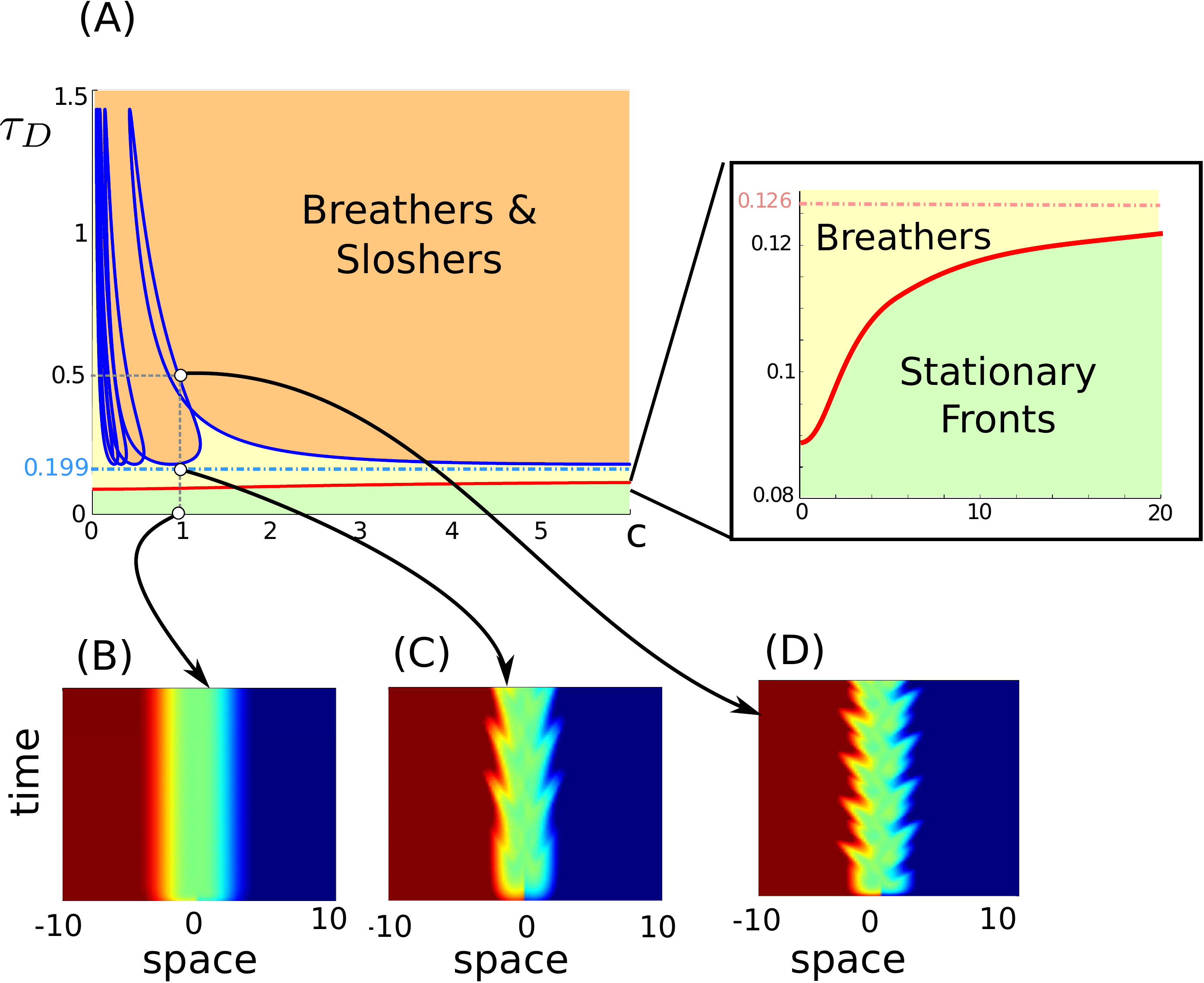}
	\caption{Stability of the front with propagation delays for an inverse Mexican hat connectivity: $w_e=1.7$, $w_i=1.2$, $\sigma_e=3$ and $\sigma_i=2$. (A) Bifurcation diagram in the $(c,\tau_D)$ plane with the first Hopf bifurcation curves given by equations \eqref{eq:SolRelation2bisProp} (red curve) and \eqref{eq:SolRelation1bis} (blue curve). These two curves asymptote the value of the delay corresponding to the purely synaptic case: $\tau_D\approx 0.126$ for \eqref{eq:SolRelation2bisProp} and $\tau_D\approx 0.199$ for \eqref{eq:SolRelation1bis}. (B) $c=1$, $\tau_D=0$: stationary front, (C): $c=1$, $\tau_D=0.2$: pulsatile front, (D) $c=1$, $\tau_D=0.5$: asymmetric complex pulsatile activity. }
	\label{fig:PropDelayFront}
\end{figure}

We represent in figure \ref{fig:PropDelayFront} the two bifurcation curves obtained from the two relationships of the above analysis: equations \eqref{eq:SolRelation2bisProp} (red curve) and \eqref{eq:SolRelation1bis} (blue curve) in the $(c,\tau_D)$ plane. We first note that the second relationship \eqref{eq:SolRelation2bisProp} always provides the first Hopf curve as it is bellow the blue curve for all values of $c$. We also remark that the first relationship~\eqref{eq:SolRelation1bis} produces a curve with several loops accumulating when $c \to 0$. These curves naturally divide the $(c,\tau_D)$ plane into three regions where the fronts can be stationary and stable or unstable with a pulsatile behavior with destabilizing modes being either symmetric (see panel (C) of Fig.~\ref{fig:PropDelayFront}) or asymmetric (see panel (D) of Fig.~\ref{fig:PropDelayFront}). Finally, we remark that constant delays are always needed in order to destabilize, through a Hopf bifurcation, non-monotonous stationary fronts such as the ones studied here.

%% file: Conclusion.tex
In this paper, we showed, using neural fields formalism, that the delays inherent to the transmission of information in the cortex can induce instabilities yielding to localized pulsatile activity. Such behaviors, related to the presence of delay-induced Hopf bifurcations, intrinsically depend on a number of parameters of the model such as the shape of the connectivity function, the threshold of the firing rate function or the strength of the input applied to the network. Moreover, these phenomena fundamentally depend on the dimension of the neural field. The cortical sheet is essentially two-dimensional, yet most studies investigate the behavior of neural fields in one dimension as a benchmark for evidencing the role of different parameters in the response of neural fields. Here, we investigated the behavior of neural fields, receiving spatially-dependent input, in the presence of delays, in an overarching formalism independent of the dimension of the space on which the neural field is posed. This allowed to uncover the dependence of observed phenomena on the dimension of the neural field. 

For localized input, we showed that the resulting pulse (or bump) solution may loose stability through a Hopf bifurcation due to the delays. In the one-dimensional neural field with inverse Mexican hat connectivity (locally inhibitory, laterally excitatory), we showed that the minimal constant delay necessary to destabilize the pulse solution tends to zero when the input strength approaches the minimal value of the input triggering a super-threshold pulse. In other words, if the input amplitude is not large enough, the resulting pulse is very sensitive to the presence of delays, and even small delays, as those involved in the synaptic transmission of information, can destabilize the formation of the pulse, yielding pulsatile activity. This phenomenon also has implications in numerical simulations: indeed, time-discretization of the neural-field equations intrinsically involve delays, and the phenomenon of destabilization of the pulse for arbitrarily small delays prevents from finding stable numerical schemes in these regions of parameters. Interestingly, this phenomenon is essentially one-dimensional: delays corresponding to the destabilization of the pulse in dimension higher than one are lowerbounded by a constant value, and pulses arising in dimension two for instance are much more stable. Moreover, we have seen that the shape of the interaction kernel shapes the form of the destabilization. Even in the class of inverse Mexican hat connectivity, different amplitudes and extents of the interaction show qualitatively distinct properties. In particular, while breathing pulses appear as the symmetric deformation mode destabilizes in most cases, asymmetric modes, corresponding to sloshing pulses, may or may not destabilize, depending on the shape of the interconnection map, and particularly on the amplitude of the connectivity kernel. 

Similar phenomena were observed for front solutions. Fronts are essentially one-dimensional phenomena: along one direction, the solution goes from a high value to a low-value. These fronts may be monotonous or not, depending on the properties of the interaction kernel. In the monotonous case, we showed that fronts are unconditionally stable in one dimension. It may be natural to think therefore that the projection of the solution on the front axis in higher dimension may also remain stable, and therefore that no destabilization may arise even in dimension higher than one. However, in higher dimensions, more phenomena may arise. In particular, an instability may arise in transverses directions of the fronts. We show that such transversal instabilities do arise, and take the example of a two-dimensional neural field model. These phenomena are not limited to fronts. Similar phenomena may arise for partially localized input, such as in the case of elongated pulses. This is a natural perspective of the present work. 

The interplay between constant delays and propagation delays was investigated. Our formalism provide the closed-form formulae for Hopf bifurcation curves as a function of the speed of the transmission as well as the value of the delay. We investigated the shape of these curves in one dimension in the case of the pulse and front. We observed interesting qualitative distinctions between the pulse and front case. Indeed, while in both case, we observe that in the limit of infinite propagation speed, the value of the Hopf bifurcation point corresponding to constant delays case is recovered, new phenomena appear at low speeds. In the example treated for the pulse in one dimension, we observe that the propagation delay acts as increasing the effective delay, and the value of constant delay necessary to destabilize the pulse reaches zero: below a certain speed, no constant delay is necessary to destabilize the pulse. More surprising is the emergence of an unstable asymmetric (sloshing) mode at low values of the speed, while in the system with no propagation delays, the asymmetric mode is unconditionally stable. Interestingly, these observations do not generalize to the case of the front. In that case, constant delays are essential to destabilize the front, and even at very small transmission speeds, the front is stable in the absence of constant delays. In the case treated, both symmetric and asymmetric modes can be destabilized by the presence of constant delays, and the Hopf bifurcation curves do asymptote to these values in the limit of infinite propagation speed. While the symmetric mode necessitates, at low speed, smaller constant delays to destabilize, the asymmetric mode yields a complex Hopf bifurcation curve (as a secondary bifurcation), and the stability of the asymmetric mode is a complex function of both constant delays and information propagation speed. 

This work therefore extends previous works on the destabilization of spatially localized activity in several directions. Most previous studies dealing with localized pattern formation of activity in neural fields dealt with one dimensional neural fields, and the destabilization occurred as a consequence of the timescale ratio between neural activity and an additional variable, typically an the adaptation variable coupled to the neuronal dynamics following the Pinto and Ermentrout model~\cite{pinto-ermentrout:01,pinto-ermentrout:01b}. Depending on the timescale of the adaptation variable, destabilization of the localized pulse may arise, both in the Mexican hat and inverse Mexican hat models~\cite{folias-bressloff:04,F}, and this additional variable can destabilize a symmetric and asymmetric mode, nonlinear analysis was developed showing that on the center manifold, the Pinto-Ermentrout model can be reduced to the normal form of the Hopf bifurcation, as was indicated by the linear analysis. It is not surprising that the Pinto-Ermentrout model and the delayed neural-field equations display qualitatively the same type of behaviors. Indeed, the adaptation variable acts as a distributed delay term on the activity variable (with an exponential delay kernel). But in contrast with the synaptic and propagation delay phenomenon, this delayed interaction corresponds systematically to a negative feedback loop whatever the interaction kernel. This is why the Pinto-Ermentrout model may display breathing and sloshing pulses even when the connectivity is locally excitatory, which is not the case of the one-dimensional delayed neural field with local excitation we analyzed here. 

Reduction to normal form around the bifurcations found is a natural extension of the present work, and would allow characterizing the stability of the oscillatory pattern beyond the Hopf bifurcation point. In contrast with~\cite{F}, the presence of delays necessitates to work in complex spaces with a specific quadratic form for the delays~\cite{hale-lunel:93}. Methods for the reduction to normal form in neural-field equations have recently been developed~\cite{veltz2013,vangils-etal:13}. While the first reference uses a $\mathbbm{L}^2$ approach, the second uses a semi-group approach (sun star formalism). This would allow for instance to uncover the presence of codimension two bifurcations, such as Bautin bifurcations, corresponding to changes in the stability of the periodic orbits generated by the Hopf bifurcations. Such bifurcations do arise for instance in the Pinto-Ermentrout model, as nicely shown by Folias in~\cite{F}. Whether such instabilities do arise in biological tissues remains to be explored. Of course, both adaptation and transmission delays arise in such systems. Moreover, stochastic effects add up. In that context, it may be reasonable to consider neural field equations arising from mean-field reductions of stochastic neuronal networks~\cite{T1,T}. The main difference is the fact that the voltage to rate transforms will no more be step functions, which will make the system much more complex to analyze, in particular because of the fact that the equilibria (stationary pulses and fronts) may no more be explicitly derivable. This study therefore provides a first step in the direction of understanding the role of constant and propagation delays in the response of neural fields to localized input, and evidences the importance of taking into account the dimension on which the system is posed, even in systems displaying an essentially one-dimensional response. 